\documentclass[english]{siamart171218}
\usepackage[latin9]{inputenc}
\usepackage{verbatim}
\usepackage{calc}
\usepackage{amsmath}
\usepackage{amssymb}
\usepackage{graphicx}
\usepackage{hyperref}
\usepackage{comment}
\usepackage{lipsum}
\usepackage{algorithm,algorithmic}
\usepackage{subfig}
\usepackage{float}
\usepackage{xcolor}

\makeatletter
\newtheorem{thm}{\protect\theoremname}
  \theoremstyle{plain}
  \newtheorem{prop}[thm]{\protect\propositionname}
  \theoremstyle{plain}
  
\theoremstyle{plain}
\def \r {\right}
\def \l {\left}

\def \I {\mathbb{I}}
\def \CC {\mathbb{C}}
\def \LL {\mathcal{L}}
\def \BB {B_+(\mathbb{I})}
\def \RR {\mathbb{R}}

\makeatother

\usepackage{babel}
  \providecommand{\lemmaname}{Lemma}
  \providecommand{\propositionname}{Proposition}
\providecommand{\theoremname}{Theorem}

\newcommand{\note}[1]{\textbf{Note: #1}}
\newcommand{\edit}[1]{{ #1}}
\newcommand{\edita}[1]{{ #1}}

\hypersetup{
   colorlinks=true,%
   citecolor=black,%
   filecolor=black,%
   linkcolor=black,%
   urlcolor=black
}

\begin{document}

\title{\edita{Sparse Inverse Problems Over Measures:\\ Equivalence of the Conditional Gradient and Exchange Methods}}
\author{Armin Eftekhari and Andrew Thompson\thanks{AE and AT have contributed equally to this work. AE is with the Institute of Electrical Engineering at the \'{E}cole Polytechnique F\'{e}d\'{e}rale de Lausanne, Switzerland. AT is with the National Physical Laboratory, United Kingdom.
}}

\maketitle
\begin{abstract}
    We study an optimization program over nonnegative Borel measures that encourages sparsity in its solution. Efficient solvers for this program are in increasing demand, as it arises when learning from data generated by a ``continuum-of-subspaces'' model, a recent trend with applications in signal processing, machine learning, and \edita{high-dimensional} statistics. We prove that the conditional gradient method (CGM) applied to this infinite-dimensional program, as proposed recently in the literature, is equivalent to the exchange method (EM) applied to its Lagrangian dual, which is a semi-infinite program. In doing so, we formally connect such infinite-dimensional programs to the well-established field of semi-infinite programming. 

On the one hand, the equivalence established in this paper allows us to provide a rate of convergence for EM which is more general than those existing in the literature. On the other hand, this connection and the resulting geometric insights might in the future lead to the design of improved variants of CGM for infinite-dimensional programs, which has been an active research topic. CGM is also known as the Frank-Wolfe algorithm. 
\end{abstract}

\section{Introduction\label{sec:Problem-Setup}}

We consider the following affinely-constrained optimization over nonnegative Borel measures:
\begin{equation}\label{eq:swapped}
\begin{cases}
\displaystyle\min_{x}&L\left(\displaystyle\int_{\I}\Phi(t)x(dt)-y\right)\\ 
\textrm{subject to}&
\|x\|_{TV} \le 1\\
& x\in\BB.
\end{cases}
\end{equation}
Here, $\I$ is a compact subset of Euclidean space,   $\BB$ denotes all nonnegative Borel measures supported on $\I$, and 
\begin{equation}
\|x\|_{TV} = \int_{\I} x(dt)
\label{eq:TV}
\end{equation}
is the \emph{total variation} of measure $x$, see for example~\cite{eftekhari2018sparse}.\footnote{\edita{It is also common to define the TV norm as half of the right-hand side of \eqref{eq:TV}, see \cite{gibbs2002choosing}.}}  We are particularly interested in the case where $L:\CC^m\rightarrow\RR$ is a differentiable \emph{loss function} and $\Phi:\I\rightarrow\CC^m$ is a continuous function.  
Note that Program~\eqref{eq:swapped} is an \emph{infinite-dimensional} problem and that the constraints ensure that the problem is bounded.  \edita{In words}, Program~\eqref{eq:swapped}  searches for a nonnegative measure on $\I$ that minimizes the loss above, while controlling its total variation. 
%In this context, Problem (\ref{eq:swapped}) searches for a measure $\hat{x}$ on $\I$ that minimizes the loss $L(\int_{\I} \Phi(t)x(dt)-y)$, while constraining its total variation. Other choices for $w(t)$ can be viewed as a weighted total variation constraint. 
This problem and its variants have  received significant attention~\cite{schiebinger2015superresolution,boyd2017alternating,candes2013super,candes2014towards,de2012exact,fernandez2013support,denoyelle2015support} in signal processing and machine learning, 
%connection with the ubiquitous problem in data and signal processing of reconstructing a sum of a few weighted sources from noisy linear measurements, with applications including superresolution imaging. 
see Section~\ref{sec:Applications} for more details.

It was recently proposed in~\cite{boyd2017alternating} to solve Program  (\ref{eq:swapped}) using the celebrated \emph{conditional gradient method} (CGM)~\cite{frank1956naval}, also known as the Frank-Wolfe algorithm, adapted to optimization over nonnegative Borel measures. The CGM algorithm minimizes a differentiable, convex function over a compact convex set, and proceeds by iteratively minimizing linearizations of the objective function over the feasible set, generating a new descent direction in each iteration. The classical algorithm performs a descent step in each new direction generated, while in the \emph{fully-corrective} CGM, the objective is minimized over the subspace spanned by all previous directions~\cite{holloway1974FC}. It is the fully-corrective version of the algorithm which we consider in this paper.

It was shown in~\cite{boyd2017alternating} that, when applied to Program \eqref{eq:swapped}, CGM generates a sequence of finitely supported measures, with a single parameter value $t^l\in\I$ being added to the support in the $l$th iteration. Moreover, \cite{boyd2017alternating} established that the convergence rate of CGM  here is $\mathcal{O}\left(\frac{1}{l}\right)$, where $l$ is the number of iterations, thereby extending the standard results for finite-dimensional CGM. A full description of CGM and its convergence guarantees can be found in Section~\ref{sec:CGM}.

On the other hand, the (Lagrangian) dual of Program \eqref{eq:swapped} is a finite-dimensional optimization problem with infinitely many constraints, often referred to as a \emph{semi-infinite program} (SIP), namely
\begin{equation}
\begin{cases}
\displaystyle\max_{\lambda,\alpha }&\text{Re}\left\langle \lambda,y\right\rangle -L_{\circ}\left(-\lambda\right)-\alpha\\
\textrm{subject to}&\text{Re}\left\langle \lambda,\Phi(t)\right\rangle \le\alpha,\qquad t\in \I\\
& \alpha \ge 0,\end{cases}\label{eq:dual of swapped general}
\end{equation}
where $\langle\cdot,\cdot\rangle$ denotes the standard Euclidean inner product over $\CC^m$. Above, 
\begin{equation}
L_{\circ}(\lambda)=\sup_{z\in\mathbb{C}^{m}}\text{Re}\langle\lambda,z\rangle-L(z)
\label{eq:cvx cnj}
\end{equation}
denotes the Fenchel conjugate of $L$. As an example, when $L(\cdot)=\frac{1}{2}\|\cdot\|_2^2$, it is easy to verify that $L_\circ = L$. 
%In this case, the objective function of \eqref{eq:dual of swapped general} becomes $\text{Re}\langle\lambda,z\rangle-\frac{\|\lambda - y\|_2^2}{2} -\alpha$. 
For the sake of completeness,  we verify the duality of Programs \eqref{eq:swapped} and \eqref{eq:dual of swapped general} in Appendix \ref{sec:infinite duality proof}. Note that the Slater's condition for the finite-dimensional Program \eqref{eq:dual of swapped general} is met and there is consequently no duality gap between the two Programs (\ref{eq:swapped}) and (\ref{eq:dual of swapped general}).

There is a large body of research on  SIPs such as Program~\eqref{eq:dual of swapped general}, see for example~\cite{hettich1993review,lopez2005review,floudas2007adaptive}, and 
%There are three main categories of methods that have become popular for solving 
%SIPs: discretization methods, exchange methods and localization methods. 
 we are particularly interested in solving Program \eqref{eq:dual of swapped general} with \emph{exchange methods}. In one instantiation -- which for ease we will refer to as \emph{the} exchange method (EM) -- one forms a sequence of nested subsets of the constraints in Program \eqref{eq:dual of swapped general}, adding in the $l$th iteration a single new constraint corresponding to the parameter value $t^l\in \I$ that maximally violates the constraints \edita{of Program \eqref{eq:dual of swapped general}}. The  finite-dimensional problem with these constraints is then solved and the process repeated. Convergence of EM has been established under somewhat general conditions, but results concerning rate of convergence are restricted to more specific SIPs, see Section~\ref{sec:EM} for a full description of the EM.

\paragraph{Contribution} The main contribution of this paper is to establish that, for Program~\eqref{eq:swapped} and provided the loss function $L$ is both strongly smooth and strongly convex, \emph{CGM and EM are dual-equivalent}. More precisely, the iterates of the two algorithms produce the same objective value and the same finite set of parameters in each iteration; for CGM, this set is the support of the current iterate of CGM and, for EM, this set is the choice of constraints in the dual program.

The EM method can also be viewed as a \emph{bundle method} for Program~\eqref{eq:dual of swapped general} as discussed in Section~\ref{sec:related work}, and the duality of CGM and bundle methods is well known for finite-dimensional problems. This paper establishes dual-equivalence in the emerging context of optimization over measures on the one hand and the well-established semi-infinite programming on the other hand. 
%In so doing, we formally connect such infinite-dimensional programs to the well-established field of semi-infinite programming. 

On the one hand, the equivalence established in this paper allows us to provide a rate of convergence for EM which is more general than those existing in the literature;   see Section~\ref{sec:related work} for a thorough discussion of the prior art.  On the other hand, this connection and the resulting geometric insights might lead to the design of improved variants to CGM, another active research topic \cite{boyd2017alternating}.

\paragraph{Outline}

We begin in Section~\ref{sec:Applications} with some motivation, describing the key role of Program \eqref{eq:swapped} in data and computational sciences. Then in Sections~\ref{sec:CGM} and~\ref{sec:EM}, we give a more technical introduction to CGM and EM, respectively. We present the main contributions of the paper in Section \ref{sec:Equivalence}, establishing the dual-equivalence of CGM and EM for Problems~\eqref{eq:swapped} and~\eqref{eq:dual of swapped general}, and deriving the rate of convergence for EM. Related work is reviewed in Section \ref{sec:related work} and \edit{some geometric insights into the inner workings of CGM and EM are provided in Section \ref{sec:geometry}. We conclude this paper with a discussion of \edita{the} future research directions.}

\section{Motivation}\label{sec:Applications}

Program \eqref{eq:swapped} has diverse applications in data and computational  sciences. In {signal processing} for example,  each $\Phi(t)\in\mathbb{C}^m$ is an \emph{atom} and the  set of all atoms $\{\Phi(t)\}_{t\in \I}$ is sometimes referred to as the \emph{dictionary}. In radar applications, for instance, $\Phi(t)$ is a copy of a known \emph{template}, arriving at time $t$. In this context, we are interested in \emph{signals} that have a \emph{sparse} representation in this dictionary, namely signals that  can be written as the superposition of a small number of atoms. Any such signal $\dot{y}\in\mathbb{C}^m$ can be written as  
\begin{equation}
\dot{y} = \int_{\I} \Phi(t) \dot{x}(dt),
\label{eq:atoms}
\end{equation}
where $\dot{x}$ is a \emph{sparse} measure, selecting the atoms that form  $\dot{y}$. More specifically, 
\begin{equation}
\dot{x} = \sum_{i=1}^k \dot{a}_i \cdot \delta_{\dot{t}_i},
\label{eq:sparse measure}
\end{equation}
for an integer $k$, positive \emph{amplitudes} $\{\dot{a}_i\}_{i=1}^k$, and \emph{parameters} $\{\dot{t}_i\}_{i=1}^k\subset \I$. 
Here, $\delta_{\dot{t}_i}$ is the Dirac measure located at $\dot{t}_i\in\I$.  We can therefore rewrite \eqref{eq:atoms} as 
\begin{equation}
\dot{y} = \int_{\I} \Phi(t) \dot{x}(dt) = \sum_{i=1}^k \Phi\l(\dot{t}_i\r) \cdot \dot{a}_i. 
\label{eq:measurements in superres nonoise}
\end{equation}
In words, $\{\dot{t}_i\}_{i}$ are the  parameters that construct the signal $\dot{y}$ and $\I$ is the \emph{parameter space}. 
We often receive $y\in\mathbb{C}^m$, a noisy copy of $\dot{y}$, and our objective in signal processing is to estimate the hidden parameters  $\{\dot{t}_i\}_i$, given the noisy copy $y$.  See Figure \ref{fig:superresExample} for an example. 

\begin{center}
\begin{minipage}{1\linewidth}
\begin{figure}[H]
\begin{center}
\subfloat[\label{fig:Picture1}]{\protect\includegraphics[width=0.49\textwidth]{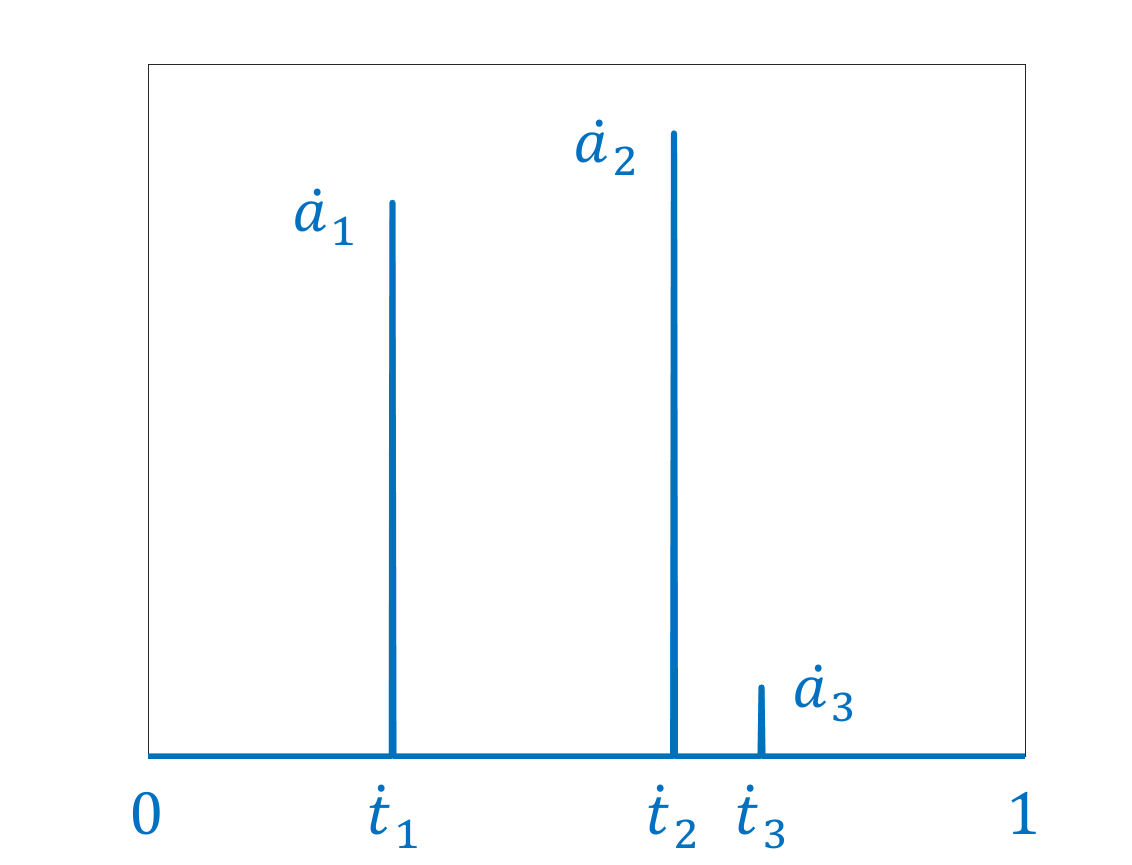}
}
\subfloat[\label{fig:Picture2}]{\protect\includegraphics[width=0.49\textwidth]{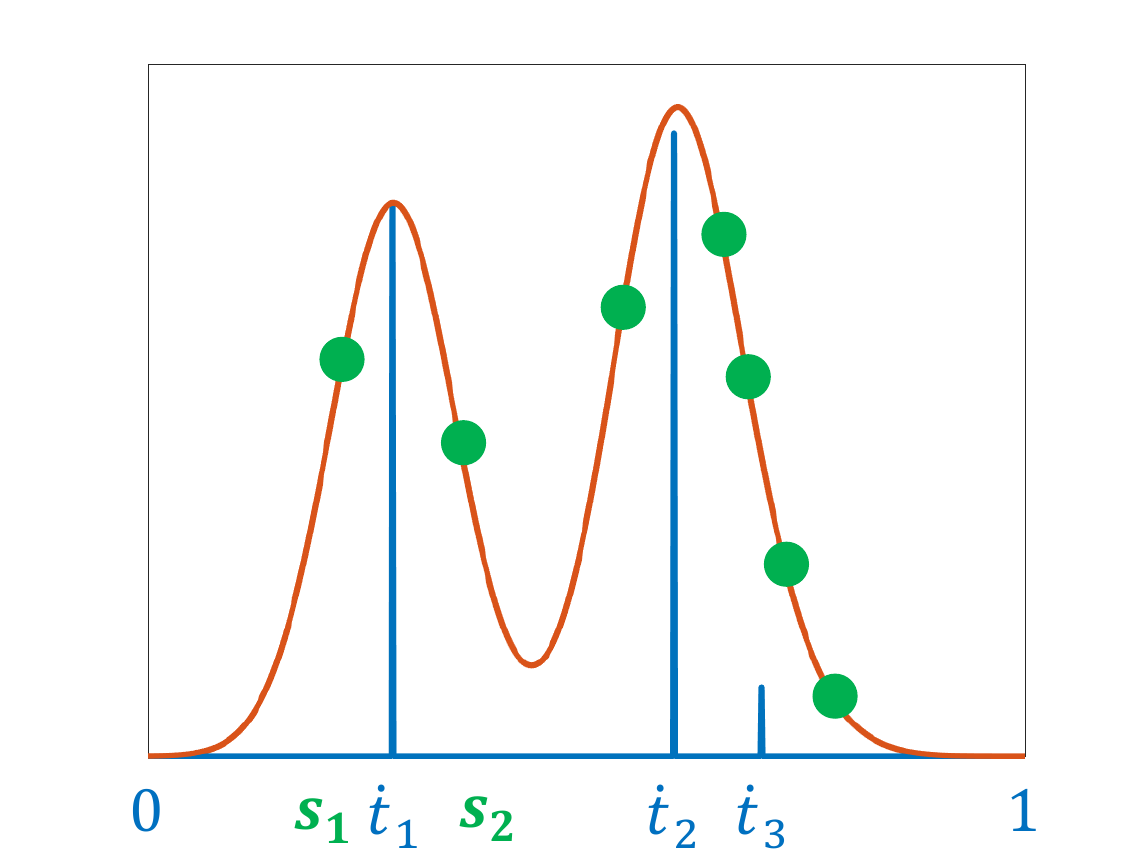}
}
\caption{\label{fig:superresExample} In this numerical example, (a) depicts the measure $\dot{x}$, see \eqref{eq:sparse measure}. Let $\phi(t)=e^{-100t^2}$ be a Gaussian window. With the choice of sampling locations $\{s_j\}_{j=1}^m\subset [0,1]$ and $\Phi(t) = [\phi(t-s_j)]_{j=1}^m \in\RR^m$, (b) depicts $\dot{y}\in\RR^m$, see \eqref{eq:measurements in superres nonoise}. Note that the entries of $\dot{y}$ are in fact samples of $(\phi\star x)(s) = \int_{\I} \phi(t-s)x(dt)$ at locations $s\in \{s_j\}_{j=1}^m$, which forms the red curve in (b). Our objective is to estimate the locations $\{\dot{t}_i\}_{i=1}^k$ from $\dot{y}$. (Given an estimate of the locations, the amplitudes $\{\dot{a}_i\}_{i=1}^k$ can also be estimated with a simple least-squares program.) This is indeed a difficult task: Even given the red curve $\phi\star \dot{x}$ (from which $\dot{y}$ is sampled), it is hard to see that there is an impulse located at $\dot{t}_3$. Solving Program \eqref{eq:swapped} with $\|x\|_{TV}\le b$ for large enough $b$ uniquely recovers $x$, as proved in \cite{eftekhari2018sparse}. In this paper, we describe Algorithms \ref{fig:eq form of CGM} and \ref{fig:EM} to solve Program \eqref{eq:swapped}, and establish their equivalence. }
\end{center}
\end{figure}
\end{minipage}
\end{center}

To that end, Program \eqref{eq:swapped} searches for a nonnegative measure $\widehat{x}$ supported on $\I$ that minimizes the {loss} $L(\int_{\I} \Phi(t)x(dt)-y)$, while \edita{encouraging its} \emph{sparsity} through the {total variation} constraint $\|x\|_{TV}\le 1$.  Under certain conditions on $\Phi$ and when $L=\frac{1}{2}\|\cdot\|_2^2$, a minimizer $\widehat{x}$ of Program \eqref{eq:swapped} is a robust estimate of the true measure $\dot{x}$ in the sense that $d(\widehat{x},\dot{x}) \le c \cdot L(y-\dot{y}) $ for a known factor $c$ and in a certain metric $d$ \cite{eftekhari2018sparse,duval2017characterization,de2012exact,Candes2014}. 

The super-resolution problem outlined above is an example of learning under a ``continuum-of-subspaces'' model, in which data belongs to the union of infinitely many subspaces. For super-resolution in particular, each subspace corresponds to fixed locations $\{t_i\}_{i=1}^K$.  This model is a natural generalization of the ``union-of-subspaces'' model, which is a central object in {compressive sensing} \cite{candes2008introduction}, wavelets \cite{mallat2008wavelet}, and feature selection in statistics \cite{hastie2015statistical}, to name a few. The use of continuum-of-subspaces models is on the rise as it potentially addresses the drawbacks of the union-of-subspaces models, see for example \cite{zhu2017approximating}. 
\edit{As another application of Program \eqref{eq:swapped}, $y$ might represent the training labels in a classification task or, in the classic moments problem, $y$ might collect the moments of an unknown distribution. Various other examples are given in \cite{boyd2017alternating}. }

Note that Program \eqref{eq:swapped} is an {infinite-dimensional} problem as the search is over all \edita{nonnegative} measures supported on $\I$. It is common in practice to restrict the support of $x$  to a uniform grid on $\I$, say $\{t_i\}_{i=1}^n\subset \I$, so that $x = \sum_{i=1}^n a_i \delta_{t_i}$ for nonnegative amplitudes $\{a_i\}_{i=1}^n$. Let $a\in \RR_+^n$ be the vector formed by the amplitudes and concatenate the vectors $\{\Phi(t_i)\}_{i=1}^n\subset \mathbb{C}^m$ to form a (usually very flat) matrix $\Phi \in\mathbb{C}^{m\times n}$. Then we may rewrite Program \eqref{eq:swapped} as 
\begin{equation}
\begin{cases}
\displaystyle\min_{a } & L\l( \Phi \cdot a - y \r)\\
\mbox{subject to} & \langle 1_n , a\rangle \le 1\\
& a\ge 0,
\end{cases}
\label{eq:lasso}
\end{equation}  
where $1_n\in\mathbb{R}^n$ is the vector of all ones. When $L(\cdot)=\frac{1}{2}\|\cdot\|_2^2$ in particular, Program \eqref{eq:lasso} reduces to the well-known \emph{nonnegative Lasso} \cite{slawski2013non}. 

\edita{The first issue with the above} 	``gridding'' approach is that there is often a mismatch between the atoms $\{\Phi(\dot{t}_i)\}_{i=1}^k$ that are present in $\dot{y}$ and the atoms listed in $\Phi$, namely $\{\Phi(t_i)\}_{i=1}^n$. As a result, $\dot{y}$ often does \emph{not} have a sufficiently sparse representation in $\Phi$. In the context of signal processing, this problem is known as the ``frequency leakage'', see Figure \ref{fig:leakage}. Countering the frequency leakage by excessively increasing the grid size $n$ leads to
\edita{increased \emph{coherence}, namely, increased similarity between the columns of $\Phi$. In turn, the statistical guarantees for finite-dimensional problems (such as Program \eqref{eq:lasso}) often deteriorate as the coherence grows~\cite[Section 1.2]{candes2011compressed}.  Loosely speaking, Program \eqref{eq:lasso} does {not} decouple the optimization error from the statistical error, and this pitfall can be avoided by directly studying the infinite-dimensional Program~\eqref{eq:swapped}, see~\cite{Candes2014}.}
\edita{Moreover, the  gridding approach  is only applicable when the parameter space $\mathbb{I}$ is low-dimensional (see the numerical example in Section~\ref{sec:CGM}), often requires post-processing \cite{tang2013sparse}, and might  lead to  numerical instability with larger grids, see Program~\eqref{eq:limited supp iteration l}. Lastly, the gridding approach ignores the continuous structure of $\mathbb{I}$ which, as discussed in Section~\ref{sec:future}, plays a key role in developing new optimization algorithms, see after \eqref{eq:main 2}. The moment technique~\cite{Candes2014,Tang2013} is an alternative to gridding for a few special choices of $\Phi$ in Program~\eqref{eq:swapped}.}

\edit{This discussion encourages us to directly study the infinite-dimensional Program \eqref{eq:swapped}; it is this direction that is pursued in this work and in \cite{boyd2017alternating,duval2017sparse,2017-Duval-IP-lasso,Eftekhari2013,support,Eftekhari2013a}. Indeed, this direct approach  provides a unified and rigorous framework, independent of gridding or its alternatives. In particular, the direct approach perfectly decouples the optimization error (caused by gridding, for instance) from the statistical error of Program~\eqref{eq:swapped}, and matches the growing trend in statistics and signal processing that aims at providing theoretical guarantees for directly learning the  underlying (continuous) parameter space $\I$~\cite{chandrasekaran2012convex,bredies2013inverse,candes2014towards,Eftekhari2013a,Tang2013}.
}

\begin{center}
\begin{minipage}{1\linewidth}
\begin{figure}[H]
\begin{center}
\subfloat[\label{fig:Picture3}]{\protect\includegraphics[width=0.33\linewidth]{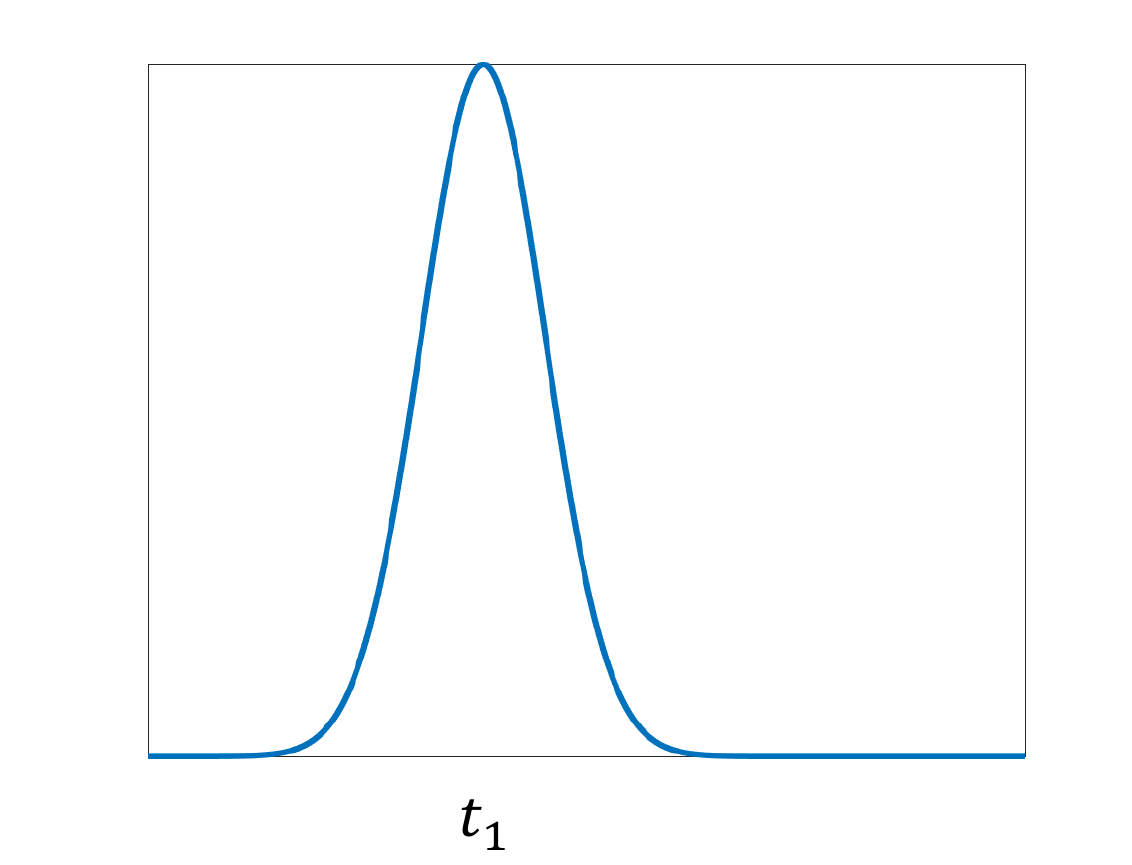}
}
\subfloat[\label{fig:Picture4}]{\protect\includegraphics[width=0.33\linewidth]{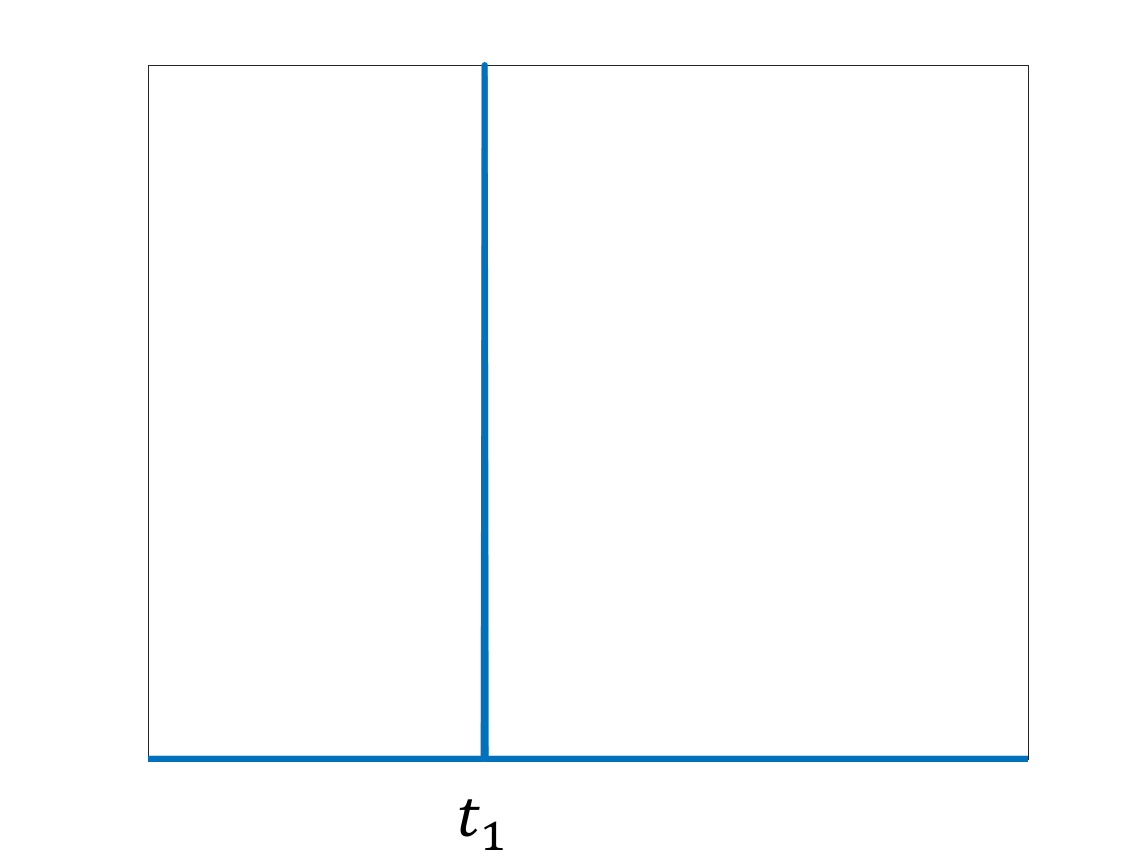}
}

\subfloat[\label{fig:Picture5}]{\protect\includegraphics[width=0.33\linewidth]{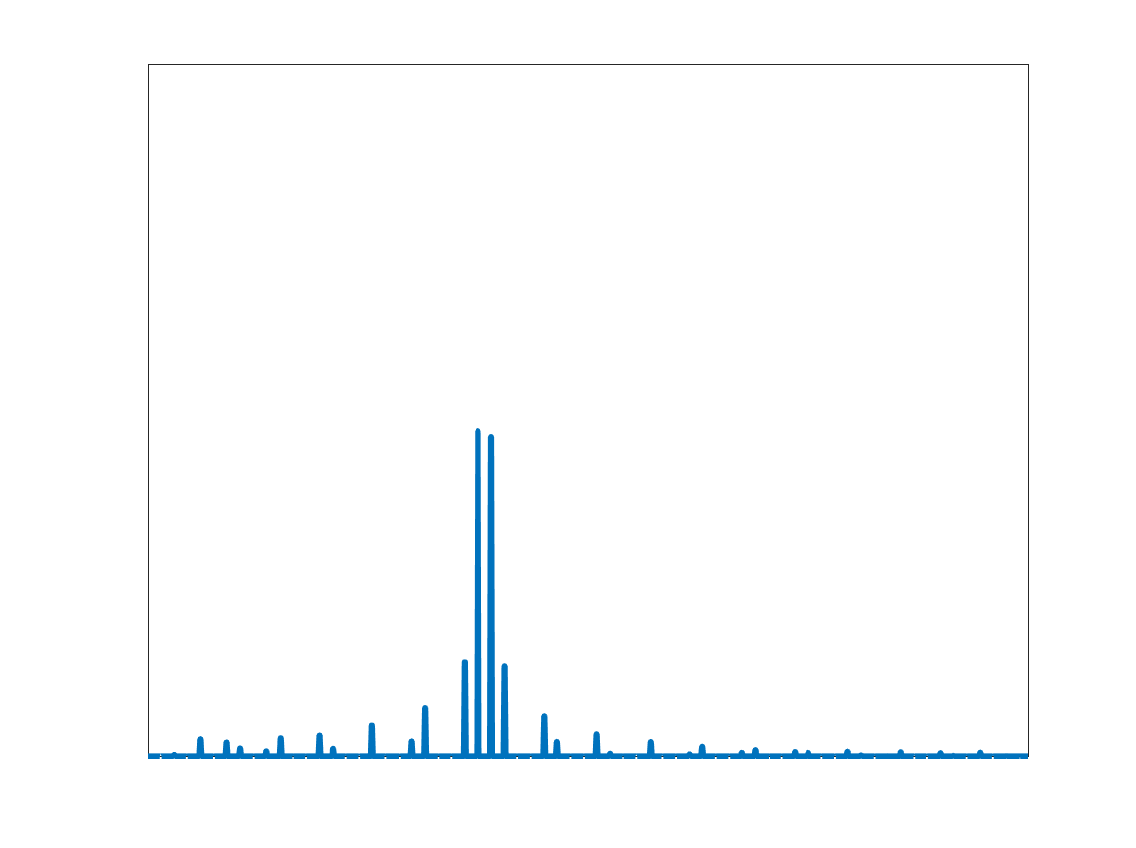}
}
\caption{\label{fig:leakage} (a) depicts a translated Gaussian window, namely, $\phi(t-t_1)=e^{-100(t-t_1)^2}$ for translation $t_1\in [0,1]$. Equivalently, $\phi(t-t_1) = (\phi\star\delta_{t_1})(t)$, as represented in (b). On the other hand, (c) shows the coefficients of the least-squares approximation of the translated window $\phi(t-t_1)$ in the dictionary $\{\phi(t-i/N)\}_{i=1}^{N}$ for $N=66$. By comparing (b) and (c), we observe that $\phi(t-t_1)$ loses its sparse representation after gridding. \edita{See the discussion at the end of Section \ref{sec:Applications} for more details.} 
%This problem is alleviated by choosing finer and finer uniform grids for the interval $[0,1]$ but at the cost of numerical instability in practice. This observation encourages directly solving the infinite-dimensional Program \eqref{eq:swapped}, as opposed to discretizing it.
}
\end{center}
\end{figure}
\end{minipage}
\end{center}

\section{Conditional Gradient Method} \label{sec:CGM}

In this section and the next one, we review two algorithms for solving Program \eqref{eq:swapped}. The first one is the conditional gradient method~\cite{frank1956naval}, a popular first-order algorithm for constrained optimization. The popularity of CGM partly stems from the fact that it is projection free, unlike projected gradient descent, for example, which requires projection onto the feasible set in every iteration. 

More specifically, CGM solves the general constrained optimization problem
$$
\displaystyle\min_{x\in \mathcal{F}} f(x) 
$$
where $f(x)$ is a differentiable function and $\mathcal{F}$ is a compact convex set. Given the current iterate $x^{l-1}$, CGM finds a search direction $s^l$ which minimizes the linearized objective function, namely, $s^l$ is a solution to
\begin{equation}
\displaystyle\min_{s\in \mathcal{F}}  f(x^{l-1})+ \langle s-x^{l-1}, \nabla f(x^{l-1} ) \rangle
\label{eq:linearized}
\end{equation}
Note that we may remove the additive terms independent of $s$ without changing the minimizers of Program \eqref{eq:linearized}.  The classical CGM algorithm then takes a step along the direction  $s^l-x^{l-1}$, namely $$x^l = x^{l-1}+\gamma^l \cdot (s^l-x^{l-1}),$$ for some step size $\gamma^l \in (0,1]$.   In a similar spirit, {fully-corrective} CGM chooses $x^l$ within the convex hull of all previous update directions \cite{holloway1974FC}. To be specific, fully-corrective CGM (which we simply refer to as CGM henceforth) sets $x^l$ to be a minimizer of 
$$\begin{cases}
\displaystyle\min &f(x) \\
\mbox{subject to} & x\in\textrm{conv}(s^1,\ldots,s^l). 
\end{cases}$$
In the context of sparse regression and classification, CGM is particularly appealing because it produces sparse iterates. Indeed, because the objective function in Program \eqref{eq:linearized} is linear in $s$, there always exist a minimizer of Program (\ref{eq:linearized}) that is an extreme point of the feasible set $\mathcal{F}$.  In our case, \edit{we have that}
\begin{equation*}
\mathcal{F} = \l\{ x\in \BB :  \|x\|_{TV} \le 1  \r\},
\end{equation*}
and any extreme point of $\mathcal{F}$ is therefore of the form $ \delta_{t}$ with $t\in \I$. It follows that each iterate $x^l$ of CGM is at most $l$-sparse, namely, supported on a subset of $\I$ of size at most $l$. 

In light of the discussion above, CGM applied to (\ref{eq:swapped}) is summarized in Algorithm \ref{fig:eq form of CGM}. Note that we might interpret Algorithm \ref{fig:eq form of CGM} as follows. Let $x_p$ be a minimizer of Program \eqref{eq:swapped}, supported on the index set $T_p \subset \I$. If an \emph{oracle} gave us the correct support $T_p$, we could  have recovered $x_p$ by solving Program~\eqref{eq:swapped} restricted to the support $T_p$ rather than $\I$. Since we do not have access to such an oracle, at iteration $l$, Algorithm \ref{fig:eq form of CGM}
\begin{enumerate}
\item finds an atom $\Phi(t^l)$ that reduces the objective of Program \eqref{eq:swapped} the most, namely an atom that is least correlated with the gradient at the current residual $\int_{\I} \Phi(\tau) x^{l-1}(d\tau) - y$, and then
\item adds  $t^l$  to the support. 
% $T^{l-1}$ to form the
%updated support $T^{l}$. 
\end{enumerate}
When $L(\cdot)=\frac{1}{2}\|\cdot\|_{2}^{2}$ in particular, Algorithm \ref{fig:eq form of CGM} reduces to the well-known \emph{orthogonal matching pursuit} (OMP) for sparse regression \cite{tropp2007signal}, adapted to measures.  

\begin{center}
\begin{algorithm}
\vspace{5pt}
\textbf{Input:} Compact set $\I$, continuous function $\Phi:\I\rightarrow\mathbb{C}^{m}$, differentiable function $L:\mathbb{C}^{m}\rightarrow\mathbb{R}$, vector $y\in\mathbb{C}^{m}$, and tolerance $\eta\ge 0$.\\ \\
\textbf{Output:} Nonnegative measure $\widehat{x}$ supported on $\I$. \\ \\
\textbf{Initialize:} Set $l=1$, $T^0=\emptyset$, and $x^{0} \equiv 0$. \\ \\
\textbf{While} $\|\nabla L ( \int_{\I} \Phi(\tau)x^{l-1}(d\tau) - y) \|_2>\eta$, \textbf{do} 
\begin{enumerate}
\item Let $t^{l}$ be a minimizer of 
\begin{equation}
\min_{t\in \I}  \l\langle\Phi(t),\nabla L\left(\int_{\I}\Phi(\tau)x^{l-1}(d\tau)-y\right)\r\rangle.
\label{eq:before simplification}
\end{equation}
\item Set $T^{l}=T^{l-1}\cup \{t^{l}\}$.
\item Let $x^{l}$ be a minimizer of 
\begin{equation}
\begin{cases}
\displaystyle\min_{x}&L\left(\displaystyle\int_{\I}\Phi(t)x(dt)-y\right)\\ 
\textrm{subject to}&\|x\|_{TV}\le 1\\
& \textrm{supp}(x)\subseteq T^l\\
& x\in \BB.
\end{cases}
\label{eq:limited supp iteration l}
\end{equation}
\end{enumerate}
\textbf{Return:} $\widehat{x}=x^l$. 
\caption{CGM for solving Program \eqref{eq:swapped}}
\label{fig:eq form of CGM}
\end{algorithm}
\end{center}

\edit{The convergence rate of CGM has been established in~\cite{boyd2017alternating}, relying heavily upon~\cite{jaggi2013revisiting}, and is reviewed next for the sake of completeness. We first note that the infinite dimensional Program~\eqref{eq:swapped} has the same optimal value  as the finite dimensional program
\begin{equation}
\displaystyle\min_{z\in C_{\I}} L(z-y),
\label{eq:finite dim2}
\end{equation}
where $C_{\I}\subset\mathbb{C}^{m}$ is the convex hull of $\{\Phi(t)\}_{t\in \I}\cup\{0\}$, namely
\begin{equation}
C_{\I} := \l\{\int_{\I} \Phi(t) x(dt) : x\in \BB, \, \|x\|_{TV} \le  1   \r\}.
\label{eq:cvxhull closure2}
\end{equation} 
Indeed, both problems  share the same objective value and their respective solutions \edita{$\widehat{z}$} and \edita{$\hat{x}$} satisfy
{$$\widehat{z}=\int_{\I}\Phi(t)\widehat{x}(dt).$$} 
It should be emphasized  that, while the problems are in this sense equivalent, solving Program~\eqref{eq:finite dim2} does not recover the underlying sparse measure but only its projection into the measurement space \edita{$\mathbb{C}^m$}. As described in Section~\ref{sec:Applications}, in many applications it is precisely the underlying sparse measure which is of interest. 
%It is for this reason that the authors of~\cite{boyd2017alternating} argue that there is additional value in considering the problem as an infinite-dimensional optimization over measures rather than simply as a finite-dimensional problem. 
A convergence result for CGM applied to Program~\ref{eq:swapped} may be obtained by first establishing that its iterates $x^i$ are related to the iterates $z^i$ of CGM applied to the finite-dimensional Program~\eqref{eq:finite dim2} by $z^l=\int_{\I}\Phi(t)x^l(dt)$. The convergence proof from~\cite{jaggi2013revisiting} can then be followed to obtain the convergence rate. Let us now turn to the details.}

For the rest of this paper, we assume that $L$ is both \emph{strongly smooth} and \emph{strongly convex}, namely, there exists $\gamma\ge 1$ such that 
\begin{equation}
\frac{\l\|x-x'\r\|_{2}^{2}}{2\gamma}\le L(x)-L(x')-\l\langle x-x',\nabla L(x')\r\rangle\le\frac{\gamma}{2}\l\|x-x'\r\|_{2}^{2},\label{eq:smooth and convex}
\end{equation}
for every $x,x'\in\mathbb{C}^{m}$.  In words, $L$ can be approximated by quadratic functions at any point of its domain. For example,  $L(\cdot)=\frac{1}{2}\|\cdot\|_2^2$ satisfies \eqref{eq:smooth and convex} with $\gamma = 1$. 
Let us also define 
\begin{equation}
r:=\max_{t\in \I}\|\Phi(t)\|_{2}.
\label{eq:def of mu}
\end{equation}
The convergence rate of Algorithm \ref{fig:eq form of CGM} is given by the following result, which is similar to the result originally given in~\cite{boyd2017alternating}, except that we replace the curvature condition in~\cite{boyd2017alternating} with the strongly smooth and convex assumption in (\ref{eq:smooth and convex}), see Appendix \ref{sec:Proof-of-Proposition cvg of CGM} for the proof. 
\begin{prop}
\textbf{\emph{(Convergence rate of Algorithm \ref{fig:eq form of CGM}) \label{prop:(Convergence-of-CMG)}}}For $\gamma\ge1$, suppose
that $L$ satisfies (\ref{eq:smooth and convex}).\footnote{Strictly speaking, strong convexity is not required for Proposition \ref{prop:(Convergence-of-CMG)}. That is, the far left term in (\ref{eq:smooth and convex})
can be replaced with zero.} Suppose that Program  (\ref{eq:before simplification}) is solved to within an accuracy of $2\gamma r^2\epsilon$ in every iteration of Algorithm \ref{fig:eq form of CGM}. Let $v_{p}$ be the optimal value of Program (\ref{eq:swapped}). Let also $v^l_{CGM}$ be the optimal value of Program (\ref{eq:limited supp iteration l}). 
Then, at iteration $l\ge1$, it holds that 
\begin{equation}
v_{CGM}^{l}-v_{p}\le\frac{4\gamma r^2(1+\epsilon)}{l+2}.
\label{eq:decay rate CGM}
\end{equation}
\end{prop}
Assuming that $L$ satisfies \eqref{eq:smooth and convex}, it is not difficult to verify that Program \eqref{eq:swapped} is a convex and strongly smooth problem. Therefore CGM achieves the same convergence rate of $1/l$ that  the projected gradient descent achieves for such problems \cite{nesterov2013introductory}. \edita{We note that, under stronger assumptions, CGM can achieves linear convergence rate \cite{garber2014faster,lacoste2015global}.}

\edit{A benefit of directly working with the infinite-dimensional Program~\eqref{eq:swapped} is that it provides a unified framework  for various finite-dimensional approximations, such as the moments method \cite{candes2014towards}.  In the context of CGM, following our discussion at the end of Section \ref{sec:Applications}, a common approach to solve Program \eqref{eq:before simplification} is to search for an $O(\epsilon)$-approximate global solution over a finite grid on $\I$, as indicated in Proposition \ref{prop:(Convergence-of-CMG)}. The tractability of this gridding approach largely depends on how smooth $\Phi(t)$ is as a function of $t$, measured by its Lipschitz constant, which we denote by  $\phi$. Roughly speaking, to find an $O(\epsilon)$-approximate global solution of Program~\eqref{eq:before simplification}, one needs to search over a uniform grid of size $O(\phi/\epsilon)^{\operatorname{dim}(\I)}$. As the dimension grows, the Lipsichtz constant $\phi$ must be smaller and smaller for this brute force search to be tractable. In some important applications, the dimension $\dim(\I)$ is in fact small. In radar, array signal processing, or imaging applications, for example, $\dim(\I)\le 2$. 

As a numerical example with $\dim(\I)=1$, let us revisit the setup in Section \ref{sec:Applications} with the choice of 
\[\dot{x} = \frac{1}{4}( \delta_{0.1\pi}+\delta_{0.2\pi}+\delta_{0.3\pi}+\delta_{0.31\pi}),\]
$$\Phi(t) = [\begin{array}{ccc}
e^{-\pi\operatorname{i}(m-1)t} & \cdots & e^{\pi \operatorname{i}  (m-1)t}
\end{array}
]^\top \in \mathbb{C}^m,$$
where $m=33$. This $\Phi$ might be considered as a  generic model for a sensing device and the resulting loss of low-frequency details~\cite{candes2014towards}. Here, $\top$ stands for vector transpose.  We solve Program \eqref{eq:swapped} by applying Algorithm \ref{fig:eq form of CGM}, where Program \eqref{eq:before simplification} therein is solved on uniform grids with sizes $\{10^2,10^3,10^4\}$. The recovery error in $1$-Wasserstein metric, namely, $d_W({x}^l,\dot{x})$, is shown in Figure \ref{fig:numerical_a}.  The same experiment is repeated in Figure \ref{fig:numerical_b} after adding additive white Gaussian noise with variance of $0.01$ to each coordinate of $\dot{y}$, see~\eqref{eq:measurements in superres nonoise}. Not surprisingly, the gains obtained from finer grids are somewhat diminished by the large noise. Both experiments were performed on a MacBook Pro (15-inch, 2017) with standard configurations. Section~\ref{sec:future} outlines a few ideas for incorporating the continuous nature of $\I$ to develop new variants of CGM that would replace the naive gridding approach above.

}
%We should perhaps note that, unlike discretizing the original Program \eqref{eq:swapped} as in Program \eqref{eq:lasso}, the above discretization approach for solving Program \eqref{eq:before simplification} is not prone to numerical instability. }
%
%\edit{From a different perspective, the discretization of Program \eqref{eq:before simplification} ensures that the computational complexity of the infinite-dimensional CGM is the same as its finite-dimensional counterpart, but without the numerical instability and the loss of sparsity caused by model mismatch, see Section~\ref{sec:Applications}.} 

\section{Exchange Method}\label{sec:EM}

EM is a well-known algorithm to solve SIPs and, in particular, Program \eqref{eq:dual of swapped general}. 
In every iteration, EM adds a new constraint out of the infinitely
many in Program (\ref{eq:dual of swapped general}), thereby forming an increasingly finer discretisation of $\I$ as the algorithm proceeds. The new constraints are added where needed most, namely, at $t\in\I$ that maximally violates the constraints in Program \eqref{eq:dual of swapped general}. In other words, a new constraint is added at $t\in\I$ that maximizes  $\text{Re}\langle\lambda^l,\Phi(t)\rangle$, where $(\lambda^l,\alpha^l)$ is the current iterate.  EM is summarized in Algorithm~\ref{fig:EM}.

Let $(\lambda_d,\alpha_d)$ be a maximizer of Program \eqref{eq:dual of swapped general}. Also assume that $T_d\subset\I$ is the set of \emph{active} constraints in Program \eqref{eq:dual of swapped general}, namely $\mbox{Re}\langle \lambda_d,\Phi(t) \rangle = \alpha_d$ for every $t\in T_d$. If an {oracle} tells us what the active constraints $T_d$ are in advance, we can simply find the optimal pair $(\lambda_d,\alpha_d)$ by solving Program \eqref{eq:dual of swapped general} with $T_d$ instead of $\I$. Alas, such an oracle is not at hand. Instead, at iteration $l$, Algorithm \ref{fig:EM}
\begin{enumerate}
\item solves Program \eqref{eq:dual of swapped general} restricted to the current  constraints $T^{l-1}$ to find $(\lambda^l,\alpha^l)$, and then
\item if  $(\lambda^l,\alpha^l)$ does not violate the constraints of Program \eqref{eq:dual of swapped general} on $\I\backslash T^{l-1}$, the algorithm terminates because it has found a maximizer of Program \eqref{eq:dual of swapped general}, namely $(\lambda^l,\alpha^l)$. Otherwise, EM adds to its support  a point $t^l\in\I$ that maximally violates the constraints of Program \eqref{eq:dual of swapped general}.
\end{enumerate}
%Section \ref{sec:T system} includes a more detailed discussion of EM that is not necessary to parse the main result of this paper. 

\begin{center}
\begin{minipage}{1\linewidth}
\begin{figure}[H]
\begin{center}
\subfloat[\label{fig:numerical_a}Noise-free]{
\protect\includegraphics[width=0.5\linewidth]{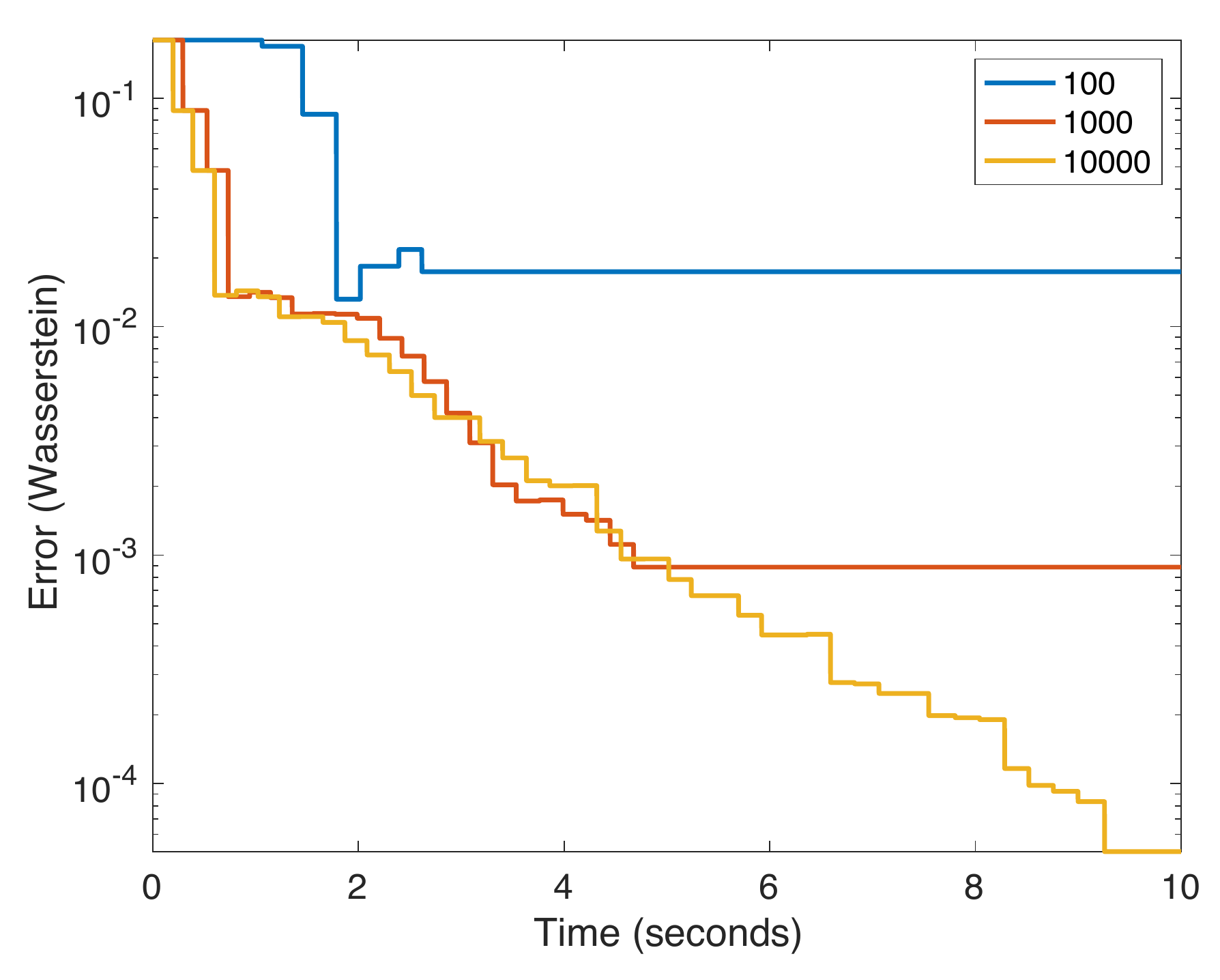}
}
\hfill
\subfloat[\label{fig:numerical_b}Noisy]{
\protect\includegraphics[width=0.5\linewidth]{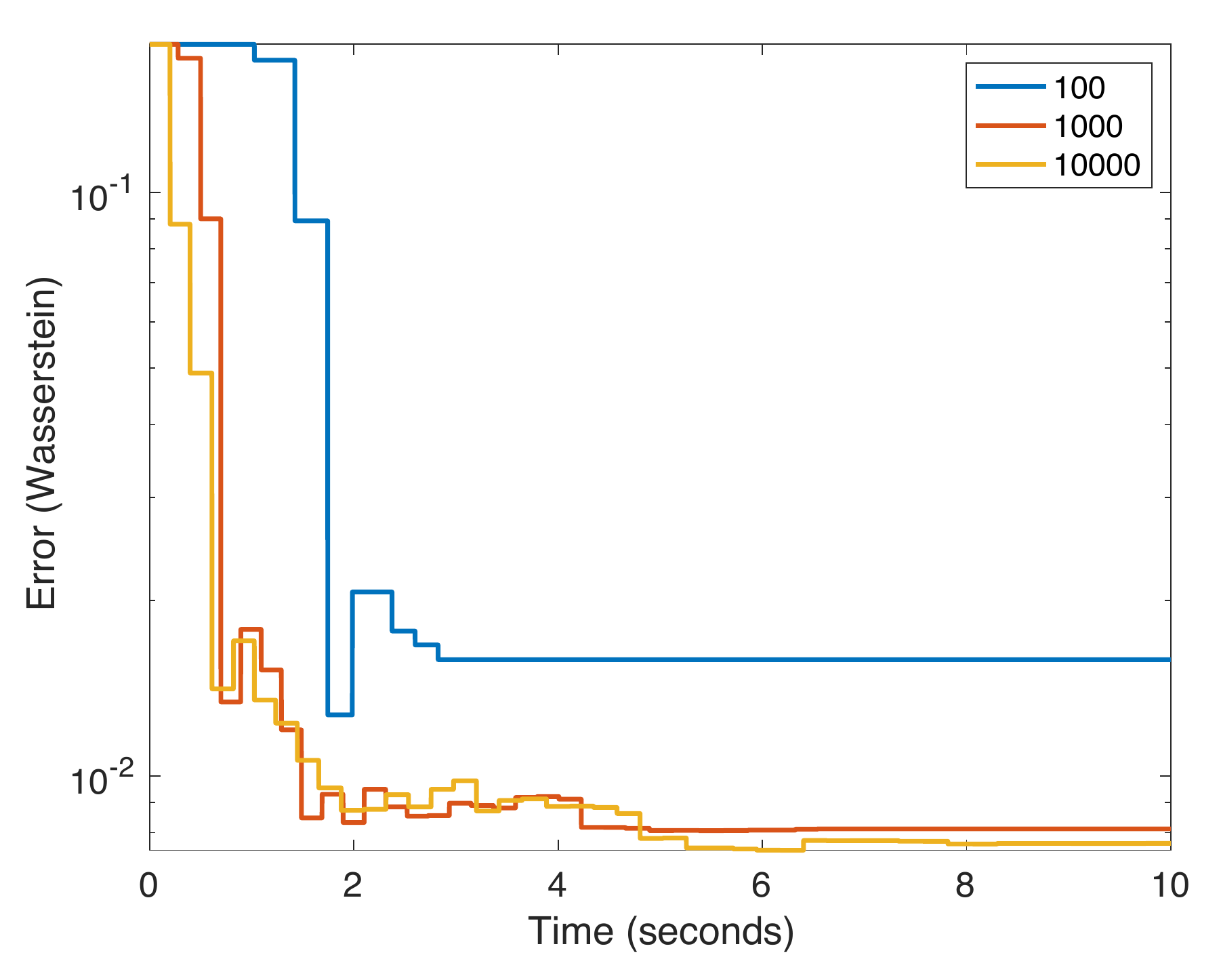}
}
\caption{Recovery error in $1$-Wasserstein metric using Algorithm \ref{fig:eq form of CGM} for the numerical example detailed at the end of Section~\ref{sec:CGM}. Grid sizes are given in the legends. \label{fig:numerical}}  
%This problem is alleviated by choosing finer and finer uniform grids for the interval $[0,1]$ but at the cost of numerical instability in practice. This observation encourages directly solving the infinite-dimensional Program \eqref{eq:swapped}, as opposed to discretizing it.
\end{center}
\end{figure}
\end{minipage}
\end{center}

\begin{algorithm}[h!]
\vspace{5pt}
\textbf{Input:} Compact set $\I$, continuous functions $\Phi:\I\rightarrow\mathbb{C}^{m}$ and
 $w:\I\rightarrow\mathbb{R}_{++}$, differentiable function $L:\mathbb{C}^{m}\rightarrow\mathbb{R}$, $y\in\mathbb{C}^{m}$ and tolerance $\eta\ge 0$.\\
 \\
\textbf{Output:} Vector $\widehat{\lambda}\in\mathbb{C}^m$ and $\widehat{\alpha}\ge 0$.\\
\\
\textbf{Initialize:} $l=1$ and $T^{0}=\emptyset$. \\
\\
\textbf{While} $\displaystyle\max_{t\in\I}\text{Re}\left\langle \lambda^{l},\Phi(t)\right\rangle>\alpha^{l}+\eta$ , \textbf{do}
\begin{enumerate}
\item Let  $(\lambda^{l},\alpha^{l})$ be a maximizer of
\begin{equation}
\begin{cases}
\displaystyle\max_{\lambda,\alpha}&\text{Re}\left\langle \lambda,y\right\rangle -L_{\circ}\left(-\lambda\right)-\alpha\\
\textrm{subject to}&\text{Re}\left\langle \lambda,\Phi(t)\right\rangle \le\alpha \qquad t\in T^{l-1}\\
& \alpha \ge 0,
\end{cases}\label{eq:main}
\end{equation}
where  $L_\circ$ is the Fenchel conjugate of $L$, see \eqref{eq:cvx cnj}.  
\item Let $t^{l} $ be the solution to
\begin{equation}\label{eq:EM maximizer}
\displaystyle\max_{t\in\I} \text{Re}\left\langle \lambda^{l},\Phi(t)\right\rangle.
\end{equation}
\item Set $T^{l}=T^{l-1}\cup t^{l}$. 
\end{enumerate}
\vspace{5pt}
\textbf{Return:}  $(\widehat{\lambda},\widehat{\alpha})=(\lambda^{l},\alpha^l)$.
\caption{EM for solving Program \eqref{eq:dual of swapped general}.}
\label{fig:EM}
\end{algorithm}

\noindent Having reviewed both CGM and EM for solving Program \eqref{eq:swapped} in the past two sections, we next establish their equivalence.

%Some preparation is required before describing the exchange method (EM) for solving Program \eqref{eq:swapped}. Consider the program 
%\begin{equation}
%\begin{cases}
%\min \,\, L\l( \int_{\I} \Phi(t) x(dt) -y \r),\\
%\int_{\I} w(t)x(dt) \le 1, 
%\label{eq:swapped limited to T}
%\end{cases}
%\end{equation}
%where the minimization is over all measures $x$ supported on a finite set $T\subset \I$. That is, Program \eqref{eq:swapped limited to T} is the restriction of Program \eqref{eq:swapped} to a finite support $T$. 

\section{Equivalence of CGM and EM} \label{sec:Equivalence}

CGM solves Program \eqref{eq:swapped} and adds a new atom in every iteration whereas EM solves the dual problem (namely Program \eqref{eq:dual of swapped general}) and  adds a new active constraint in every iteration, and both algorithms do so ``greedily''.  
Their connection goes deeper: Consider Program \eqref{eq:swapped} restricted to a finite support $T\subset \I$, namely, the program 
\begin{equation}
\begin{cases}
\displaystyle\min_{x} & L\l( \int_{\I} \Phi(t) x(dt) -y \r)\\
\mbox{subject to} & \|x\|_{TV} \le 1\\
& \mbox{supp}(x) \subseteq T\\
& x\in \BB.
\label{eq:swapped limited to T}
\end{cases}
\end{equation}
The dual of Program (\ref{eq:swapped limited to T}) is 
\begin{equation}
\begin{cases}
\displaystyle\max_{\lambda,\alpha}& \text{Re}\left\langle \lambda,y\right\rangle -L_{\circ}\left(-\lambda\right)-\alpha\\
\mbox{subject to} & \text{Re}\left\langle \lambda,\Phi(t)\right\rangle \le\alpha \qquad t\in T\\
& \alpha \ge 0.
\end{cases}\label{eq:dual of swapped}
\end{equation}
 Indeed, Program \eqref{eq:dual of swapped} is the restriction of Program \eqref{eq:dual of swapped general} to $T$. 
 Note that the complementary slackness forces any minimizer of Program \eqref{eq:swapped limited to T} to be supported on the set of active constraints of Program \eqref{eq:dual of swapped}.  Note also that Programs \eqref{eq:swapped limited to T} and \eqref{eq:dual of swapped} appear respectively in CGM and EM but with different support sets. 
The following result states that CGM and EM are in fact equivalent algorithms to solve Program~\eqref{eq:swapped}, see Appendix \ref{sec:Proof-of-Proposition equl} for the proof.
\begin{prop}
\label{prop:CMG-and-EM are equil}\textbf{\emph{(Equivalence of Algorithms \ref{fig:eq form of CGM} and \ref{fig:EM})}}  
For $\gamma\ge 1$, suppose  that $L$ satisfies (\ref{eq:smooth and convex}). Assume also that CGM and EM update their supports according to the same rule, e.g., selecting the smallest solutions if $\I\subset\mathbb{R}$. Then CGM and EM are equivalent in the sense that $T_{CGM}^{l}=T_{EM}^{l}$
for every iteration $l\ge 0$. Here, $T_{CGM}^{l}$ and $T_{EM}^{l}$ (both subsets
of $\I$) are the support sets of CGM and EM at iteration $l$, respectively. 

Furthermore, $v_{CGM}^{l}=v_{EM}^{l+1}$, where $v_{CGM}^{l}$ and
$v_{EM}^{l}$ denote the optimal values of Programs (\ref{eq:limited supp iteration l})
and (\ref{eq:main}) in CGM and EM, respectively. 
\end{prop}
The above equivalence allows us to carry convergence results from one algorithm to
another. In particular, the convergence rate of CGM in Proposition \ref{prop:(Convergence-of-CMG)} determines the convergence rate of EM, as the following result indicates, see Appendix \ref{sec:Proof-of-Proposition cvg of EM}  for the proof. 
\begin{prop}
\label{prop:cvg of EM} \textbf{\emph{(Convergence of Algorithm \ref{fig:EM})}} For $\gamma \ge 1$, suppose that $L$ satisfies (\ref{eq:smooth and convex}). Recall the definition of $r$ in (\ref{eq:def of mu}) and, for $\epsilon\ge 0$, suppose  that Program (\ref{eq:EM maximizer}) is solved to within an accuracy of $2\gamma r^2\epsilon$ in every iteration.   Let $v_{d}$ be the optimal value of Program (\ref{eq:dual of swapped general}) and $(\lambda_d,\alpha_d)$ be its unique maximizer. Likewise, let $v^l_{EM}$ be the optimal value of Program (\ref{eq:main}). 
%and recall that the iterate $(\lambda^l,\alpha^l)$ of Algorithm \ref{fig:EM} denotes its unique maximizer. 
At iteration $l\ge1$, it then holds that
\begin{equation}
v_{EM}^{l}-v_{d}\le\frac{4\gamma r^2(1+\epsilon)}{l+2},
\label{eq:cvg rate of EM}
\end{equation}
\begin{equation*}
\| \lambda^l - \lambda_d \|_2\le \sqrt{\frac{8\gamma^2 r^2(1+\epsilon)}{l+2}},
\end{equation*}
\begin{equation}
| \alpha^l - \alpha_d | \le \sqrt{\frac{8\gamma^2 r^4(1+\epsilon)}{l+2}}.
\label{eq:cvg rate of EM iterates}
\end{equation}
Furthermore,  it holds that 
\begin{align}
\max_{t\in \I} \langle \lambda^l,\Phi(t) \rangle \le \alpha_d + \sqrt{\frac{8\gamma^2 r^4(1+\epsilon)}{l+2}}.
\label{eq:feasibility}
\end{align} 
That is, the iterates $\{\lambda^l\}_l$ of Algorithm \ref{fig:EM} gradually become feasible for Program (\ref{eq:dual of swapped general}).
\end{prop}
Proposition \ref{prop:cvg of EM} states that  Program \eqref{eq:dual of swapped general}, which has  infinitely many constraints, can be solved as fast as a smooth convex program with finitely many constraints. More specifically, it is not difficult to verify that the objective function of Program (\ref{eq:dual of swapped general}) is convex and strongly smooth, see Section \ref{sec:geometry}.  Then,  \eqref{eq:cvg rate of EM} states that EM solves Program \eqref{eq:dual of swapped general} at the rate of $1/l$, the same rate at which the projected gradient descent solves a finite-dimensional problem under the assumptions of convexity and strong smoothness~\cite{nesterov2013introductory}. This is perhaps remarkable given that Program  \eqref{eq:dual of swapped general} has infinitely many constraints. Note however that the convergence of the iterates $\{(\lambda^l,\alpha^l)\}_l$ of EM to the unique maximizer $(\lambda_d,\alpha_d)$ of Program \eqref{eq:dual of swapped general} is much slower as given in \eqref{eq:cvg rate of EM iterates}, namely,  at the rate of $1/\sqrt{l}$.

%Recall that $(\lambda^l,\alpha^l)$ is the unique maximizer of Program \eqref{eq:main}. Loosely speaking, \eqref{eq:feasibility} states that $(\lambda^l,\alpha^l)$  rapidly becomes feasible for Program \eqref{eq:dual of swapped general}. 
%Indeed, the feasibility of $(\lambda^l,\alpha^l)$ for Program \eqref{eq:dual of swapped general} would also imply the optimality of this pair for Program \eqref{eq:dual of swapped general}, thereby terminating EM.
We remark that Proposition \ref{prop:cvg of EM} is novel in providing a rate of convergence for EM for a general class of nonlinear SIPs, whereas the literature on SIP only gives rates of convergence for specific problems. See Section~\ref{sec:related work} for a more detailed literature review. 

\section{Related Work} \label{sec:related work}

The conditional gradient method (CGM), also known as the Frank-Wolfe algorithm, is one of the earliest algorithms for constrained optimization~\cite{frank1956naval}. The version of the algorithm considered in this paper is the fully-corrective Frank-Wolfe algorithm, also known as the \emph{simplicial decomposition} algorithm, in which the objective is optimized over the convex hull of all previous atoms~\cite{holloway1974FC,bertsekas2011unifying}. The algorithm was proposed for optimization over measures, the context considered in this paper, in~\cite{boyd2017alternating}.

Semi-infinite programs (SIPs) have been much studied, both theoretically in terms of optimality conditions and duality theory, and algorithmically in terms of design and analysis of numerical methods for their solution. We refer the reader to the review articles~\cite{hettich1993review} and~\cite{lopez2005review} for further background. 

Exchange methods are one of the three families of popular methods for the numerical solution of SIPs, with the other two being discretisation methods and localization methods. In discretisation methods, the \edita{infinitely many} constraints are replaced by a finite subset thereof and the resulting finite dimensional problem is solved as an approximation of the SIP. In localization methods, a sequence of local (usually quadratic) approximations to the problem are solved. 

Global convergence of discretisation methods has been proved for linear SIPs~\cite{goberna1998book}, but no general convergence result exists for nonlinear SIPs~\cite{hettich1993review}. Global convergence of exchange methods has been proved for general SIPs~\cite{hettich1993review}, but to the authors' best knowledge there is no general proof of rate of convergence, except for more specific problems. For localization methods, local superlinear convergence has been proved assuming strong sufficient second-order optimality conditions, which do not hold for all SIPs~\cite{fontecilla1987localization}. The guarantees extend to global convergence of more sophisticated algorithms which combine localization methods with global search, see~\cite[Section 7.3]{lopez2005review} and references therein. We refer the reader to~\cite{hettich1993review,lopez2005review} for more details on existing convergence analysis of SIPs. {Against this background, the convergence rate of the EM, established here in Proposition \ref{prop:cvg of EM} for a wide class of nonlinear SIPs, represents a new contribution.}

\edit{The exchange method described in this paper can also be viewed as the cutting plane method, also known as Kelley's method~\cite{kelley1960cutting,Kelley1999} applied to Program~\eqref{eq:dual of swapped general}. Dual equivalence of conditional gradient methods and cutting plane methods is well known for finite-dimensional problems, see for example~\cite{bach2013equivalence,bertsekas2011unifying,lacoste-julien2013block}, and these results agree with the dual equivalence established in this paper. 

EM may also be viewed as a bundle method for unconstrained optimization \cite{bagirov2014bundle,bach2013submodular}. Bundle methods construct piecewise linear approximations to an objective function using a ``bundle'' of subgradients from previous iterations. As a special case, given a convex and smooth function $u$ and convex (but not necessarily smooth) function $v$, the function $u+v$ may be minimized by constructing piecewise linear approximations to $v$, generating the sequence of iterates $\{\lambda^l\}_l$ specified as} 
\begin{equation}
\label{sub_iter}
\lambda^l\in \arg\min_{\lambda} \l( u(\lambda)+\max_{1\le i\le l-1} \text{Re}\langle \lambda,\partial v(\lambda^i)\rangle \r),
\end{equation}
where $\partial v(\lambda^i)$ is a subgradient of $v$ at $\lambda^i$. 
To establish the connection with EM, note that Program (\ref{eq:dual of swapped general}) can be rewritten as the unconstrained problem
\begin{equation}
\max_{\lambda\in\CC^m}\,\, \text{Re}\left\langle \lambda,y\right\rangle -L_{\circ}\left(-\lambda\right)-\max_{t\in\I}\langle \lambda,\Phi(t)\rangle.
\label{eq:sup formulation}
\end{equation}
Setting $u(\lambda)=-\text{Re}\langle \lambda,y\rangle+L_\circ(-\lambda)$ and $v(\lambda)=\max_{t\in \I} \langle \lambda,\Phi(t)\rangle$, and then applying the bundle method produces the iterates
\begin{equation}
\lambda^l \in \arg\max \text{Re}\left\langle \lambda,y\right\rangle -L_{\circ}\left(-\lambda\right)-\max_{t\in T^l}\langle \lambda,\Phi(t)\rangle \equiv \mbox{Program \eqref{eq:main}}.
\end{equation}
That is, EM applied to Program \eqref{eq:dual of swapped general} and the bundle method described above applied to Program~\eqref{eq:sup formulation} produce the same iterates. The dual equivalence of CGM and the bundle method has previously been noted for various finite dimensional problems, see for example~\cite[Chapter 7]{bach2013submodular}. 
%An alternative approach to proving Proposition~\ref{prop:CMG-and-EM are equil} would be to invoke the known dual equivalence between CGM and Kelley's method for the finite dimensional problems (\ref{rewritten2}) and (\ref{eq:sup formulation}) respectively, and then to use the equivalence between Kelley's method and the EM algorithm. 
However, we are not aware of any extension of these finite-dimensional results to SIPs and their dual problem of optimization over Borel measures. {In this sense, the equivalence established in Proposition \ref{prop:CMG-and-EM are equil} is novel.}

\section{Geometric Insights} \label{sec:geometry}

\edit{This section collects a number of useful insights about CGM/EM and Program \eqref{eq:swapped} in the special case where $\I\subset\mathbb{R}$ is a compact subset of the real line and the function $\Phi:\I\rightarrow\mathbb{C}^m$ is a \emph{Chebyshev system} \cite{karlin1966tchebycheff}, see Section \ref{sec:Problem-Setup}.} 
\begin{definition}\label{def:Tsys}\textbf{\emph{(Chebyshev system)}} Consider a compact interval $\I\subset\mathbb{R}$ and a continuous function $\Phi:\I \rightarrow\mathbb{C}^m$. Then $\Phi$ is a Chebyshev system if $\{\Phi(t_i)\}_{i=1}^m \subset \mathbb{C}^m$ are linearly independent vectors for any choice of distinct $\{t_i\}_{i=1}^m\subset \I$.\footnote{Note that Definition \ref{def:Tsys} is slightly different from the standard one in \cite{karlin1966tchebycheff} which requires $\Phi$ to be real-valued.} 
\end{definition}
Chebyshev systems are widely used in classical approximation theory and generalize the notion of ordinary polynomials. Many functions form Chebyshev systems, for example sinusoids or translated copies of the Gaussian window, and we refer the interested reader to \cite{karlin1966tchebycheff,MarkovBook,eftekhari2018sparse} for more on their properties and applications. 
 Let $C_{\I}\subset\mathbb{C}^{m}$ be the convex hull
of $\{\Phi(t)\}_{t\in \I}\cup\{0\}$, namely 
\begin{equation}
C_{\I} := \l\{\int_{\I} \Phi(t) x(dt) : x\in \BB, \, \|x\|_{TV} \le  1   \r\}.
\label{eq:cvxhull closure}
\end{equation} 
Note that $\{x\in \BB : \| x\|_{TV} \le  1\}$ is a compact set.  Then, by the continuity of $\Phi$ and with an application of the dominated convergence theorem, it follows that $C_{\I}$  is a compact set too. Since $\Phi$ is by assumption a Chebyshev system, $C_{\I}\subset \CC^m$ is in fact a \emph{convex body}, namely a compact convex set with non-empty interior. 
% Since $\Phi$ is a Chebyshev system by assumption, $C_{\I}$ is in fact a \emph{convex body}, namely a compact convex subset  of $\mathbb{C}^m$ with nonempty interior.  
Introducing $z=\int_{\I}\Phi(t)x(dt)$, we note that Program (\ref{eq:swapped}) is equivalent to the program 
\begin{equation}
\displaystyle\min_{z\in C_{\I}} L(z-y).
\label{eq:finite dim}
\end{equation}
The compactness of $C_{\I}$ and the strong convexity of $L$ in \eqref{eq:smooth and convex} together imply that Program \eqref{eq:finite dim} has a unique minimizer $y_p\in C_{\I}$, which can be written as $y_p = \int_{\I} \Phi(t)x_p(dt)$, where $x_p$ itself is a minimizer of Program \eqref{eq:swapped}. 
For example, when $L(\cdot)=\frac{1}{2}\|\cdot\|_{2}^{2}$, Program (\ref{eq:finite dim}) projects $y$ onto $C_{\I}$. That is, $y_p$ is the orthogonal projection of $y$ onto the convex set $C_{\I}$.  

Given the equivalence of Programs \eqref{eq:swapped} and \eqref{eq:finite dim}, we might say that solving Program \eqref{eq:swapped}  ``denoises'' the signal $y$ from a signal processing viewpoint, in the sense that it finds a nearby signal $y_p = \int_{\I} \Phi(t) {x}_p(dt)$ that  has a sparse representation in the dictionary $\{\Phi(t)\}_{t\in \I}$ (because $x_p$ is  a sparse measure). To be more specific, by Carath\'eodory's theorem \cite{barvinok2002course}, every $y_p\in C_{\I}$ can be written as a convex combination of at most $m$  atoms of the dictionary $\{\Phi(t)\}_{t\in \I}$. 
On the other hand, the Chebyshev assumption on $\Phi$ implies that  $\{\Phi(t)\}_{t\in \I}$ are the \emph{extreme points} of $C_{\I}$~\cite[Chapter II]{karlin1966tchebycheff}. 
Here, an extreme point of $C_{\I}$ is a point in $C_{\I}$ that cannot be written as a convex combination of other points in $C_{\I}$. It then follows that this atomic decomposition of $y_p$ is unique, and $x_p$ is necessarily $m$-sparse. 
We may note the analogous result in the finite-dimensional case. Indeed, the Lasso problem is known to have a unique solution whose sparsity is no greater than the rank of the measurement matrix, provided the columns of the measurement matrix are in general position~\cite{tibshirani2012uniqueness}.

%Indeed, by Caratheodory theorem, every point in $C_{\I}$ and particularly $y_p$ can be written as the convex combination of at most $m$ points 

At iteration $l$ of  CGM, let $C^l\subset \mathbb{C}^m$ be the convex hull of $\{\Phi(t)\}_{t\in T^{l}}\cup\{0\}$, namely 
\begin{equation}
C^l := \l\{ \sum_{t\in T^l} \Phi(t) \cdot a_t : \sum_{t\in T^l} a_t \le 1,\, a_t \ge 0,\,\, \forall t\in T^l \r\}.
\label{eq:def of Cl}
\end{equation}
Similar to the argument above, we observe that
%, namely 
%\begin{equation}
%C^l  = \left\{ \sum_{t \in T^l }a_{t}\cdot \frac{\Phi(t)}{w(t)}\,:\,\sum_{t \in T^l} a_{t}\le1,\,\,a_{t}\ge0,\,\,\forall t \in T^l\right\} \subset\mathbb{C}^{m}.
%\label{eq:def of Cl}
%\end{equation}
% Let
%also $v_{CGM}^{l}$ be the optimal value of Program (\ref{eq:limited supp iteration l}).
%Then, after introducing $z=\int_{\I}\Phi(t)x(dt)$, we write that
Program \eqref{eq:limited supp iteration l} is equivalent to
\begin{equation}
\displaystyle\min_{z\in C^l}L(z-y).
\label{eq:projection at iteration l}
\end{equation}
%In particular, if $L(\cdot)=\frac{1}{2}\| \cdot \|_2^2$, then Program \eqref{eq:projection at iteration l} projects $y$ onto the convex set $C^l$.  
As with Program \eqref{eq:finite dim}, Program \eqref{eq:projection at iteration l} has a unique minimizer $y^l\in C^l$ that satisfies $y^l = \int_{T^l} \Phi(t) x^l(dt)$ and $x^l$ is a minimizer of Program \eqref{eq:limited supp iteration l}. By Carath\'eodory's theorem again, $x^l$ is at most $m$-sparse. In other words, there always exists an $m$-sparse minimizer $x^l$ to Program \eqref{eq:limited supp iteration l}; iterates of CGM are always sparse and so are the iterates of EM by their equivalence in Proposition \ref{prop:CMG-and-EM are equil}.
%That is, a minimizer $x^l$ exists that is supported on a subset of $T^l$ of size at most $m$.  

%For  arbitrary but distinct $\{t_i\}_{i=1}^m\subset\I$, the vectors $\{\Phi(t_i)\}_{i=1}^m$ are linearly independent vectors, because $\Phi$ is by assumption a Chebyshev system. Therefore, $C_{\I}$ is a convex body, namely a compact convex set with nonempty interior. Moreover, $\{\Phi(t)/w(t)\}_{t\in \I}$ are the \emph{extreme points} of $C_{\I}$. Here, an extreem point of $C_{\I}$ cannot be written as a convex combination of other points in $C_{\I}$. By Caratheodory theorem, every $z\in C_{\I}$ can be written as a convex combination of at most $m$ extreme points of $C_{\I}$.  
%
In addition, note that the chain $C^1 \subseteq C^2 \subseteq \cdots \subseteq C_{\I}$ provides a sequence of increasingly better approximations to $C_{\I}$. CGM eventually terminates when $y_p = y^l \in C^l \subseteq C_{\I}$, which happens as soon as $C^l$ contains the \emph{face} of $C_{\I}$ to which $y^p$ belongs. It is however common to use different stopping criteria to terminate CGM when $y^l$ is sufficiently close to $y_p$.

Let us now rewrite Program \eqref{eq:dual of swapped general} in a similar way. First let $C_{\I,\circ}\subset\mathbb{C}^m$ be the \emph{polar} of $C_{\I}$, namely 
$$C_{\I,\circ}=\left\{ \lambda:\text{Re}\left\langle \lambda,z\right\rangle \le1,\,\forall z\in C_{\I}\right\} =\left\{ \lambda:\text{Re}\left\langle \lambda,\Phi(t)\right\rangle \le 1,\,\forall t\in \I\right\},$$
where the second identity follows from the definition of $C_{\I}$. Let also $g_{C_{\I,\circ}}=\gamma_{C_{\I}}$ denote the \emph{gauge function}  associated
with $C_{\I,\circ}$ and the \emph{support function} associated with $C_{\I}$, respectively \cite{rockafellar2015convex}.  That is, for $\lambda\in\mathbb{C}^m$, we define 
\begin{equation}
\label{eq:gauge}
g_{C_{\I,\circ}}(\lambda)
:=\begin{cases}
\displaystyle\min_{\alpha} &\alpha\\
\textrm{subject to}&\lambda\in\alpha\cdot C_{\I,\circ}\\
& \alpha \ge 0
\end{cases}
=\max_{t\in \I}\,\, {\mbox{Re}\langle \lambda,\Phi(t)\rangle } 
=\max_{z\in C_{\I}}\,\, {\mbox{Re}\langle \lambda,z\rangle } 
=: \gamma_{C_{\I}}(\lambda)
,
\end{equation}
where $\alpha \cdot C_{\I,\circ}=\{\alpha \lambda : \lambda \in C_{\I,\circ} \}$.
%%\footnote{The naming is not very strict since  $\|\cdot\|_{C_{\I,\circ}}$ is only a norm when $C_{\I,\circ}$ and in turn $C_{\I}$ are \emph{centrally symmetric}, namely when every $\pm \Phi(t)$  belong to our dictionary. Formally, when $\pm\Phi(t)\in\{\Phi(\tau)\}_{\tau\in \I}$ for every $t\in\I$.} 
In words, $g_{C_{\I,o}}(\lambda)=\gamma_{C_{\I}}(\lambda)$ is the smallest $\alpha$ at which the inflated ``ball'' $\alpha \cdot C_{\I,\circ}$ first reaches $\lambda$. 
By usual convention, the optimal value above is set to infinity when
the problem is infeasible, namely when the ray that passes through
$\lambda$ does not intersect $C_{\I,\circ}$. It is also not difficult to verify that $g_{C_{\I,\circ}}=\gamma_{C_{\I}}$ is a positively-homogeneous convex function.  Using \eqref{eq:gauge}, we may rewrite Program~(\ref{eq:dual of swapped general}) as
\begin{equation}
\begin{cases}
\displaystyle\max_{\lambda,\alpha}&L_{\circ}\left(-\lambda\right)+\text{Re}\left\langle \lambda,y\right\rangle -\alpha\\
\textrm{subject to}&\lambda\in\alpha\cdot C_{\I,\circ}\\
& \alpha \ge 0
\end{cases}
\equiv
\max_{\lambda\in\mathbb{C}^{m}}\,\, \text{Re}\left\langle \lambda,y\right\rangle -L_{\circ}\left(-\lambda\right)-g_{C_{\I,\circ}}(\lambda).
\label{eq:equiv form of main dual}
\end{equation}
By the assumption in \eqref{eq:smooth and convex}, $L$ is strongly smooth and consequenly  $L_\circ$ is strongly convex \cite{nesterov2013introductory}. Therefore Program \eqref{eq:dual of swapped general} has a unique maximizer, which we denote by $(\lambda_d,\alpha_d)$.  The optimality of $(\lambda_d,\alpha_d)$ also immediately implies that 
\begin{equation}
\alpha_d = g_{C_{\I,\circ}} (\lambda_d).
\label{eq:alpha d n lambda d}
\end{equation}
Thanks to Proposition \ref{prop:CMG-and-EM are equil}, we likewise define the polar of $C^{l-1}$ and the corresponding gauge function to rewrite the main step of EM in Algorithm \ref{fig:EM}, namely 
\begin{equation}\label{polar_version}
\begin{cases}
\displaystyle\max_{\lambda,\alpha}&L_{\circ}\left(-\lambda\right)+\text{Re}\left\langle \lambda,y\right\rangle -\alpha\\
\textrm{subject to}&\lambda\in\alpha\cdot C^{l-1}_{\circ}\\
& \lambda\ge 0
\end{cases}
\equiv
\max_{\lambda\in\mathbb{C}^{m}}\,\, \text{Re}\left\langle \lambda,y\right\rangle -L_{\circ}\left(-\lambda\right)-g_{C^{l-1}_{\circ}}(\lambda)
\end{equation}
and the unique minimizer of the above three programs is $(\lambda^l,\alpha^l)$, where the uniqueness again comes from the strong convexity of $L_\circ$. Similar to \eqref{eq:alpha d n lambda d}, the optimality of $(\lambda^l,\alpha^l)$ immediately implies that
\begin{equation}
\alpha^l = g_{C^l_\circ} (\lambda^l).
\label{eq:alpha l n lambda l}
\end{equation}
It is not difficult to verify that 
\begin{equation}
C^1_\circ\supseteq C^2_\circ \supseteq \cdots \supseteq C_{\I,\circ}.
\label{eq:nested duals}
\end{equation}
 That is, as EM progresses, $C^l$ gradually ``zooms into'' $C_{\I,\circ}$. As with CGM, EM eventually terminates as soon as $C_{\circ}^l$ includes the face of $C_{\I,\circ}$ to which $\lambda_d/\alpha_d$ belongs, at which point $(\lambda^l,\alpha^l)=(\lambda_d,\alpha_d)$. In light of the argument in Appendix \ref{sec:Proof-of-Proposition equl}, in every iteration, we also have that 
\begin{equation}
\l\langle y^l, \lambda^l/\alpha^l\r\rangle = 1,
\end{equation}
namely the pair $(y^l,\lambda^l/\alpha^l)\in C^l\times C^l_{\circ}$ satisfies the generalized Holder inequality $g_{C^l} \cdot g_{C^l_\circ}\le 1  $ with equality~\cite{rockafellar2015convex}. Here, $g_{C^l}$ and $g_{C^l_{\circ}}$  are the gauge functions of $C^l$ and $C^l_\circ$, respectively. 
It is worth pointing out that, with the choice of $L(\cdot)=L_{\circ}(\cdot)=\frac{1}{2}\|\cdot \|_{2}^{2}$, the maximizer of (\ref{polar_version}) is the same as the (unique) minimizer of 
\[
\min_{\lambda\in\mathbb{C}^{m}}\frac{\|\lambda-y\|_{2}^{2}}{2}+g_{C^{l-1}_{\circ}}(\lambda),
\]
which might be interpreted as a generalization of Lasso  and other standard tools for sparse denoising~\cite{hastie2015statistical}. That is,  each iteration of CGM/EM can be interpreted as a simple denoising procedure.

\edit{
\section{Future Directions} \label{sec:future}

Even though, by Proposition \ref{prop:(Convergence-of-CMG)}, CGM reduces  the objective function of Program \eqref{eq:swapped} at the rate of $O(1/l)$,  the behavior of $\{x^l\}_l$, namely, the sequence of measures generated by Algorithm \ref{fig:eq form of CGM},  is often less than satisfactory. Indeed, in practice, the greedy nature of CGM leads to adding clusters of spikes to the support, many of which are spurious. 

In this sense, all applications reviewed in Section \ref{sec:Applications} will benefit from improving the performance of CGM. In particular, a variant suggested in \cite{boyd2017alternating} follows each step of CGM with a heuristic local search. Intuitively, this modification makes the algorithm less greedy and helps avoid the clustering of spikes, described above. 

The equivalence of CGM and EM discussed in this paper might offer a new and perhaps less heuristic approach to improving CGM. From the perspective of EM, a natural improvement to Algorithm \ref{fig:EM} (and consequently Algorithm \ref{fig:eq form of CGM}) might be  obtained by replacing Program \eqref{eq:main} with 
\begin{equation}
\begin{cases}
\displaystyle\max_{\lambda,\alpha}&\text{Re}\left\langle \lambda,y\right\rangle -L_{\circ}\left(-\lambda\right)-\alpha\\
\textrm{subject to}&\text{Re}\left\langle \lambda,\Phi(t)\right\rangle \le\alpha \qquad t\in T_\delta^{l-1}\\
& \alpha \ge 0,
\end{cases}\label{eq:main 2}
\end{equation}
where $T^{l-1}_\delta \subseteq \I $ is the $\delta$-neighborhood of the current support $T^{l-1}$, namely, all the points in $\I$ that are within $\delta$ distance of the set $T^{l-1}$.

At the first glance, Program \eqref{eq:main 2} is itself a semi-infinite program and not  any easier to solve than Program \eqref{eq:dual of swapped general}. However, if $\delta$ is sufficiently small  compared to the Lipschitz constant of $\Phi$, then one might use a local approximation for $\Phi$ to approximate Program \eqref{eq:main 2} with a finite-dimensional problem.  For instance, if $\Phi$ is differentiable, one could use the first order Taylor expansion of $\Phi$ around each impulse in $T^{l-1}$. As another example, suppose that $\Phi:\I=[0,1)\rightarrow\mathbb{C}^m$ and specified as 
%
%suppose $\Phi$ is differentiable and let $D\Phi(t)$ be its Jacobian at $t\in \I$. Then, for sufficiently small $\delta$, one can instead consider the program 
% \begin{equation}
%\begin{cases}
%\displaystyle\max_{\lambda,\alpha}&\text{Re}\left\langle \lambda,y\right\rangle -L_{\circ}\left(-\lambda\right)-\alpha\\
%\textrm{subject to}& \text{Re}\langle\lambda,\Phi(t) \rangle + \text{Re}\left\langle \lambda,D\Phi(t) (t'-t)\right\rangle \le\alpha \qquad  \|t'-t\| \le \delta, \,\, t\in T^{l-1}\\
%& \alpha \ge 0.
%\end{cases}\label{eq:main 3}
%\end{equation}
%It is now possible to further simplify Program \eqref{eq:main 3} by removing any dependence on $t'$. If $\|\cdot\|_*$ denotes the dual norm of $\|\cdot\|$, Program \eqref{eq:main 3} is equivalent to 
%\begin{equation}
%\begin{cases}
%\displaystyle\max_{\lambda,\alpha}&\text{Re}\left\langle \lambda,y\right\rangle -L_{\circ}\left(-\lambda\right)-\alpha\\
%\textrm{subject to}& \text{Re}\langle\lambda,\Phi(t) \rangle + \|(D\Phi(t))^{\operatorname{H}}\lambda\|_* \delta \le\alpha \qquad   t\in T^{l-1}\\
%& \alpha \ge 0,
%\end{cases}\label{eq:main 3}
%\end{equation}
%where $(D\Phi(t))^{\operatorname{H}}$ stands for the complex conjugate of the Jacobian at $t$. 
%There are also other ways to locally approximate $\Phi$. Consider, for example, the case where 
\begin{equation}
\Phi(t) = [
\begin{array}{ccc}
e^{-\pi\operatorname{i}(m-1)t} & \cdots & e^{\pi \operatorname{i}  (m-1)t}
\end{array}
]^\top \in \mathbb{C}^m,
\end{equation}
see the numerical test at the end of Section~\ref{sec:CGM}.  It is easy to verify that $\{\Phi(j/m)\}_{j=0}^{m-1}$ form an orthonormal basis for $\mathbb{C}^m$. Even though we may represent $\Phi(t)$ in Program \eqref{eq:main 2} within this basis for any $t\in \I$, this representation is not ``local'' \cite{Eftekhari2013a,zhu2017approximating}. A better local representation of $\Phi(t)$ within a $\delta$-neighborhood is obtained through the machinery of \emph{discrete prolate spheroidal wave functions} \cite{Slepian1978}.

In light of this discussion, an interesting future research direction might be to study variants of Program \eqref{eq:main 2} and their potential impact in various applications.
}

\section*{Acknowledgements}

The authors are grateful to Mike Wakin, Jared Tanner, Mark Davenport, Greg Ongie, St\'ephane Chr\'etien, and Martin Jaggi for their helpful feedback and comments. The authors are grateful to the anonymous referees for their valuable suggestions. \edita{For this project, AE was} supported by the Alan
Turing Institute under the EPSRC grant EP/N510129/1 and also by the Turing Seed Funding grant SF019.

\bibliographystyle{unsrt}
\bibliography{References}

\appendix

\section{Proof of Proposition \ref{prop:(Convergence-of-CMG)}}\label{sec:Proof-of-Proposition cvg of CGM}

\edit{
Recall the equivalent form of Program \eqref{eq:swapped} given by Program \eqref{eq:finite dim}, and let $v_p$ be the optimal value of both these programs. We first establish that the iterates $x^i$ of CGM applied to \eqref{eq:swapped} are related to the iterates $z^i$ of CGM applied to~\eqref{eq:finite dim2} by $z^i=\int_{\I}\Phi(t)x^i(dt)$. In this regard, suppose that $z^i=\int_{\I}\Phi(t)x^i(dt)$ and let $t^{i+1}$ be the solution to Program~\eqref{eq:before simplification}. Let $s^i$ be the output of the linear minimization step of CGM applied to Program~\eqref{eq:finite dim}. Then
$$s^i=\arg\min_{s\in\mathbb{C}_{\I}}\langle s,\nabla L(z^i-y)\rangle=\Phi(t^{i+1}),$$
which shows that the linear steps of CGM for both formulations coincide.}
%\begin{comment}
%Similarly, recall $C^l\subset\CC^m$, defined in \eqref{eq:def of Cl}, and the equivalent form of Program \eqref{eq:swapped limited to T} in Algorithm \ref{fig:EM} https://www.overleaf.com/project/given by Program \eqref{eq:projection at iteration l}, letting $v^l_{CGM}$ denote the optimal value of both of these programs. 
%\end{comment}
%Convergence of CGM can be established by noting that Program~\eqref{eq:swapped} is equivalent to the finite-dimensional program
%\begin{equation}\label{rewritten2}
%v_p = 
%\displaystyle\min_{z\in C_{\I}}L(z-y).
%\end{equation}
%Likewise, let $C^l \subseteq C_{\I}$ be the convex hull of $\{\Phi(t)\}_{t\in T^l}\cup \{0\}$, namely 
%\begin{equation*}
%C^l = \l\{ \sum_{t\in T^l} \Phi(t)\cdot a_t\,:\, \sum_{t\in T^l} a_t \le 1, \,\,a_t \ge 0,\,\, \forall t\in T^l \r\}.
%\end{equation*}
%Program \eqref{eq:before simplification}  is then equivalent to the finite-dimensional program 
%\begin{equation}\label{rewritten3}
%v^l_{CGM} = 
%\displaystyle\min_{z\in C^l} L(z-y)
%\end{equation}
%Consider the finite-dimensional formulation, namely Program (\ref{rewritten2}), 

\edit{Now suppose that Program (\ref{eq:before simplification}) is solved to an accuracy of $\theta \cdot \epsilon$ in every iteration, where}
\begin{equation}
\theta=
\begin{cases}
\displaystyle\sup_{\rho,z,s}& 
\frac{2}{\rho^2}\left(L(z'-y)-L(z-y)-\langle z'-z,\nabla L(z-y)\rangle\right)\\
& z' = z+\rho(s-z)\\
& z,s\in C_{\I}\\
& \rho\in [0,1].
\end{cases}
\label{eq:theta}
\end{equation}
Then we may invoke Theorem 1 in \cite{jaggi2013revisiting} to obtain that
$$v_{CGM}^{l}-v_{p}\le\frac{2\theta(1+\epsilon)}{l+2}.$$
Let us next bound $\theta$ in terms of the known quantities. Due to the assumption (\ref{eq:smooth and convex}), for any feasible pair $(z,z')$ in \eqref{eq:theta}, we have that 
\begin{align}\label{strong_bound}
L(z'-y)-L(z-y)-\langle z'-z,\nabla L(z-y)\rangle& \le\frac{\gamma}{2}\|z'-z\|_2^2 \nonumber\\
& \le \frac{\gamma\rho^2}{2} \|s-z\|_2^2
\qquad \l(z' = z +\rho(s-z) \r) \nonumber\\
& \le \frac{\gamma}{2} ( \|s\|_2^2 + \|z\|_2^2) \nonumber\\
& \le \gamma \rho^2 r^2,
\qquad \mbox{(see \eqref{eq:theta})}
\end{align}
which immediately implies that $\theta\le 2\gamma r^2 $,
%Meanwhile, given $z,s\in C_{\I}$, we have
%\begin{equation}\label{cauchy_bound}
%\|z-s\|_2^2\le 2(\|x\|_2^2+\|s\|_2^2)\le 2 r^2.
%\qquad \mbox{(see \eqref{eq:def of mu})}
%\end{equation}
%For any $\rho\in [0,1]$ and $z'=z+\rho(s-z)$, we may therefore combine (\ref{strong_bound}) and (\ref{cauchy_bound}) to obtain that
%$$f(z)-f(x)-\langle z-x,\nabla f(x)\rangle\le\gamma\rho^2 r^2.$$
%We may therefore bound $\mathcal{C}$ as
%$$\mathcal{C}\le\frac{2}{\rho^2}\cdot\gamma\rho^2 r^2=2\gamma r^2$$
and (\ref{eq:decay rate CGM}) now follows.

\section{Duality of Programs (\ref{eq:swapped}) and (\ref{eq:dual of swapped general})} \label{sec:infinite duality proof}

We show here that the dual of Program (\ref{eq:swapped}) is Program (\ref{eq:dual of swapped general}). We first observe that Program (\ref{eq:swapped}) is equivalent to 
\begin{equation}\label{rewritten}
\begin{cases}
\displaystyle\min_{z,x}&L(z-y)\\
\textrm{subject to}&z=\displaystyle\int_{\I}\Phi(t)x(dt)\\
&\|x\|_{TV}\le 1\\
& x\in \BB.
\end{cases}
\end{equation}
Introducing Lagrange multipliers $\lambda\in\CC^m$ and $\alpha\geq 0$ for the two respective constraints, the Lagrangian $\LL(x,\lambda,\alpha)$ for Program (\ref{rewritten}) is
$$\begin{array}{rcl}\LL(x,\lambda,\alpha)&=&L(z-y)-\textrm{Re}\left\langle\lambda,\left(z-\displaystyle\int_{\I}\Phi(t)x(dt)\right)\right\rangle+\alpha\left(\|x\|_{TV}-1\right)\\
&=&L(z-y)+\textrm{Re}\langle\lambda,z\rangle+\displaystyle\int_{\I}\left(\alpha -\textrm{Re}\langle\lambda,\Phi(t)\rangle\right)x(dt),\end{array}$$
and so the dual of Program (\ref{rewritten}) is
$$\max_{\lambda \in\CC^m,\alpha\geq 0}\left\{\inf_{z\in\CC^m}\left[L(z-y)+\textrm{Re}\langle\lambda,z-y\rangle\right]+\inf_{\mu\in\BB}\left[\int_{\I}\left(\alpha -\textrm{Re}\langle\lambda,\Phi(t)\rangle\right)x(dt)\right]+\text{Re}\langle\lambda,y\rangle-\alpha\right\}$$
where $\BB$ is the set of all nonnegative Borel measures supported on $\I$.
Using the definition of the Fenchel conjugate in (\ref{eq:cvx cnj}), the above problem is equivalent to 
$$\begin{cases}
\displaystyle\max_{\lambda,\alpha}&-L_{0}(-\lambda)+\textrm{Re}\langle\lambda,y\rangle-\alpha\\
\textrm{subject to}&\textrm{Re}\langle\lambda,\Phi(t)\rangle\le\alpha \qquad t\in\I\\
& \alpha \ge 0,
\end{cases}$$
which is Program (\ref{eq:dual of swapped general}).

\section{Proof of Proposition \ref{prop:CMG-and-EM are equil}} \label{sec:Proof-of-Proposition equl}

By construction, $T_{CGM}^0 = T_{EM}^0 = \emptyset$. 
Fix iteration $l\ge1$ and assume that $T_{CGM}^{l-1}=T_{EM}^{l-1}=T^{l-1}$.
We next show that $T_{CGM}^{l}=T_{EM}^{l}=T^l = T^{l-1}\cup \{ t^l\}$,
namely the two algorithms add the same point $t^l$ to their support sets
in iteration $l$. We opt for a geometric argument here that relies heavily on Section \ref{sec:geometry}. 

Recall that Program \eqref{eq:limited supp iteration l} is equivalent to Program \eqref{eq:projection at iteration l}. Recall also that $y^l = \int_{T^l} \Phi(t) x^l(dt)$ is the unique minimizer of Program \eqref{eq:projection at iteration l}, where  $x^l$ is a minimizer of Program \eqref{eq:limited supp iteration l}. 
%Let $C^{l-1} \subset \mathbb{C}^m$ be the convex hull of $\{\Phi(t)\}_{t\in T^{l-1}} \cup \{0\}$.  From Appendix \ref{sec:Duality-of-Programs}, recall that~\eqref{eq:limited supp iteration l} is equivalent to
%\begin{equation}
%\begin{cases}
%\displaystyle\min_{z\in\CC^m}&L(z-y)\\
%\textrm{subject to}&z\in C^{l-1},
%\end{cases}
%\label{eq:dual of EM prg at iteration l}
%\end{equation}
%under the change of variables $z=\int_{\I} \Phi(t)x(dt)$. Recall also that $x^l$ denotes a minimizer of \eqref{eq:limited supp iteration l} and therefore $y^l = \int_{\I} \Phi(t) x^l(dt)$ is a minimizer of \eqref{eq:dual of EM prg at iteration l}. By assumption, $L$ is strongly convex, see \eqref{eq:smooth and convex}. Therefore $y^l$ is in fact the unique minimizer of \eqref{eq:dual of EM prg at iteration l}. 
On the other hand, recall that Program \eqref{eq:main} is equivalent to Program \eqref{polar_version}, and both programs have the unique minimizer $(\lambda^l,\alpha^l)$. Since Program \eqref{eq:main} only has linear constraints, Slater's condition is met and 
%
%Let also $C^{l-1}_\circ\subset \mathbb{C}^m$ be the polar of $C^{l-1}$, namely $C^{l-1}_\circ = C_{T^{l-1},\circ}$ as defined in \eqref{eq:def of polar}.  From Appendix \ref{sec:Duality-of-Programs}, recall that~\eqref{eq:main} is equivalent to
%\begin{equation}
%\begin{cases}
%\displaystyle\max_{\lambda\in\CC^m,\alpha\geq 0}&L_{\circ}\left(-\lambda\right)+\text{Re}\left\langle \lambda,y\right\rangle -\alpha\\
%\textrm{subject to}&\lambda\in\alpha\cdot C_{T,\circ}.
%\end{cases}
%\label{eq:equiv form of EM iteration l proof}
%\end{equation}
%Recall that  $(\lambda^l,\alpha^l)$ denotes a minimizer of \eqref{eq:main}. By assumption in \eqref{eq:smooth and convex}, $L$ is strongly smooth. Therefore $L_\circ$ is strongly convex and $(\lambda^l,\alpha^l)$ is the unique maximizer of \eqref{eq:equiv form of EM iteration l proof} \cite{nesterov2013introductory}.   Moreover, since~\eqref{eq:dual of EM prg at iteration l} is a convex problem with affine constraints, 
there is no duality gap between Programs (\ref{eq:projection at iteration l}) and (\ref{polar_version}). Furthermore,
%  
% 
% 
% The following simple result, proved in Appendix \ref{sec:proof of technical lem}, is about the duality gap.  
%\begin{lem}\label{lem:technical lem}
%There is no duality gap between Programs (\ref{eq:limited supp iteration l}) and (\ref{eq:main}) if $\Phi$ is a Chebyshev system. Moreover,  $y^l=\int_{\I}\Phi(t) x^l(dt) $ is the unique minimizer of Program (\ref{eq:limited supp iteration l}), and $(\lambda^l,\alpha^l)$ is the unique maximizer of Program  (\ref{eq:main}). 
%\end{lem}
%
%
%From Appendix \ref{sec:Duality-of-Programs}, also recall that, as shown in Appendix \ref{sec:Duality-of-Programs}, Programs
%(\ref{eq:limited supp iteration l}) and (\ref{eq:main}).  
%
%
%In that appendix, we also saw that Programs \eqref{eq:limited supp iteration l}  and  \eqref{eq:main} are equivalent to Programs \eqref{eq:atomic version 2} and \eqref{eq:duality in polar}, respectively, with $T=T^{l-1}_{CGM}=T^{l-1}_{EM}$. It is more convenient to work with the latter programs in this proof.  
%
%
%
%By assumption, $L$ is strongly
%convex and therefore $y^{l}=:\int_{\I}\Phi(t)x^{l}(dt) \in\mathbb{C}^m$ is the
%unique solution of Program (\ref{eq:atomic version 2}) with $T=T^{l-1}_{CGM}=T^{l-1}_{EM}$. 
%On the other hand, the assumption that $L$ is strongly smooth implies
%that $L_{\circ}$ is strongly convex. Therefore, $(\lambda^{l},\alpha^{l})$
%is the unique solution of Program (\ref{eq:duality in polar}) with $T=T^{l-1}_{CGM}=T^{l-1}_{EM}$. 
the tuple $(y^{l},\lambda^{l},\alpha^{l})$  satisfies the KKT conditions, namely 
\[
y^l \in C^{l-1}, 
\qquad 
\lambda^l \in \alpha^l \cdot  C^{l-1}_\circ,
\qquad 
\alpha^l \ge 0,
\]
\[
\lambda^{l}=-\nabla L(y^{l}-y),\qquad \langle y^{l},\lambda^{l} \rangle=\alpha^{l},
%\qquad 
%\alpha^l y^l = 0.
\]
From the above expression for $\lambda^l$, it follows  immediately that the same point is added to the support in both Programs (\ref{eq:limited supp iteration l}) and (\ref{eq:main}), which implies that $T_{CGM}^{l}={T}_{EM}^{l}={T}^{l-1}\cup \{t^l \}$.
Finally, the above argument reveals that $v_{CGM}^{l}=v_{EM}^{l+1}$, which completes the proof of Proposition \ref{prop:CMG-and-EM are equil}.

\section{Proof of Proposition \ref{prop:cvg of EM}} \label{sec:Proof-of-Proposition cvg of EM}

%For $l\ge1$ and by the strong duality in every iteration between
%Programs (\ref{eq:limited supp iteration l}) and (\ref{eq:main}),
%we have that 
%\[
%L\left(r^{l-1}\right)=v_{CGM}^{l-1}=v_{EM}^{l}=\alpha^{l}\text{Re}\l\langle\lambda^{l},y\r\rangle-L_{\circ}\left(-\alpha^{l}\lambda^{l}\right)-\alpha^{l}.
%\]
Note that 
\begin{align}
v_{EM}^{l}-v_{d} & =v_{CGM}^{l-1}-v_{d}
\qquad \mbox{(see Proposition \ref{prop:CMG-and-EM are equil})}
\nonumber \\
 & =v_{CGM}^{l-1}-v_{p}
\qquad \mbox{(strong duality between Programs \eqref{eq:swapped} and \eqref{eq:dual of swapped general})} 
 \nonumber \\
 & \le\frac{4\gamma r^2(1+\epsilon)}{l+2},\qquad\text{(see Proposition \ref{prop:(Convergence-of-CMG)})}\label{eq:cvg of EM 1}
\end{align}
which proves the first claim in Proposition \ref{prop:cvg of EM}. 
To prove the second claim there, first recall the setup in Section \ref{sec:geometry}. Let us first show that the minimizer of Program \eqref{eq:main}, namely $(\lambda^l,\alpha^l)$, converge to the minimizer of Program \eqref{eq:dual of swapped general}, namely $(\lambda_d,\alpha_d)$. 
To that end, recall the equivalent formulation of Programs (\ref{eq:dual of swapped general},\ref{eq:main}) given in (\ref{eq:equiv form of main dual},\ref{polar_version}), and let 
\begin{equation*}
h_{C_{\I,o}}(\lambda) :=  \text{Re}\langle \lambda,y\rangle - L_\circ(-\lambda) - g_{C_{\I,o}}(\lambda),
\end{equation*}
\begin{equation}
h_{C^l_{o}}(\lambda) :=  \text{Re}\langle \lambda,y\rangle - L_\circ(-\lambda) - g_{C^{l-1}_{o}}(\lambda),
\label{eq:def of h}
\end{equation}
denote their objective functions, respectively. In particular, note that 
\begin{equation}
h_{C_{\I,\circ}}(\lambda_d)=v_d,
\qquad 
h_{C^l_{\circ}}(\lambda^l)=v^l_{EM}.
\label{eq:optimality of lambda_d}
\end{equation}
By assumption in \eqref{eq:smooth and convex}, $L$ is $\gamma$-strongly smooth and therefore $L_\circ$ is $(\gamma^{-1})$-strongly convex \cite{nesterov2013introductory}. Consequently, $-h_{C_{\I,\circ}}$ is also $(\gamma^{-1})$-strongly convex, which in turn implies that 
\begin{align}
\frac{1}{2\gamma} \| \lambda^l - \lambda_d \|_2^2 & \le -h_{C^l_{\circ}}(\lambda_d) + h_{C^l_{\circ}}(\lambda^l) + \langle \lambda_d-\lambda^l,\nabla h_{C^l_{\circ}}(\lambda^l) \rangle \nonumber\\
& = -h_{C^l_{\circ}}(\lambda_d)+ v_{EM}^l,
\qquad \mbox{(see \eqref{eq:optimality of lambda_d})}
\label{eq:cvg rate of dual pre}
\end{align}
\edit{where the inner product above disappears by optimality of $\lambda^l$ in Program \eqref{polar_version}.} 
Let us next control $h_{C^l_{\circ}}(\lambda_d)$ in the last line above by noting that 
\begin{align}
h_{C^l_{\circ}}(\lambda^d) & = 
\text{Re}\langle \lambda_d,y\rangle - L_\circ(-\lambda_d) - g_{C^l_{o}}(\lambda_d) 
\qquad \mbox{(see \eqref{eq:def of h})}
\nonumber\\
& \ge \text{Re}\langle \lambda_d,y\rangle - L_\circ(-\lambda_d) - g_{C_{\I,o}}(\lambda_d)
\qquad \l( C^l_\circ \supseteq C_{\I,\circ} \mbox{ in \eqref{eq:nested duals}}\r) \nonumber\\
& = h_{C_{\I,\circ}}(\lambda_d) \nonumber\\
& =  v_d. \qquad \mbox{(see \eqref{eq:optimality of lambda_d})}
\end{align} 
By substituting the bound above back into \eqref{eq:cvg rate of dual pre}, we find that 
\begin{align}
\| \lambda^l - \lambda_d \|_2^2&  \le 2\gamma (v^l_{EM}-v_d ) \nonumber\\
& \le \frac{8\gamma^2 r^2(1+\epsilon)}{l+2}.
\qquad \mbox{(see \eqref{eq:cvg of EM 1})}
\label{eq:cvg rate of lambdas}
\end{align}
The above bound also allows us to find the convergence rate of $\alpha^l$ to $\alpha_d$. Indeed, note that 
\begin{align}
\label{eq:cvg of alphas proof}
| \alpha^l - \alpha_d | & = \l| g_{{C^l_\circ}}(\lambda^l) - g_{C_{\I,\circ}}(\lambda_d) \r|
\qquad \mbox{(see (\ref{eq:alpha l n lambda l},\ref{eq:alpha d n lambda d}))}
\nonumber\\
& = \l| \max_{t\in T^l} \langle \lambda^l,\Phi(t)  \rangle - \max_{t\in \I} \langle \lambda_d,\Phi(t)  \rangle \r| 
\qquad \mbox{(see \eqref{eq:gauge})}
\nonumber\\
& \le \max_{t\in \I} \l| \langle \lambda^l-\lambda_d, \Phi(t)  \rangle \r| \nonumber\\
& \le \| \lambda^l - \lambda_d \|_2 \max_{t\in \I} \| \Phi(t)\|_2 \nonumber\\
& \le \sqrt{ \frac{8\gamma^2 r^2 (1+\epsilon)}{l+2}} \cdot r.
\qquad \mbox{(see (\ref{eq:cvg rate of lambdas},\ref{eq:def of mu}))}
\end{align}
With an argument similar to \eqref{eq:cvg of alphas proof}, we also find that 
\begin{align}
\l| \max_{t\in \I} \langle \lambda^l,\Phi(t)  \rangle - \alpha_d \rangle \r| 
 \le \sqrt{ \frac{8\gamma^2 r^2 (1+\epsilon)}{l+2}} \cdot r,
\end{align}
which completes the proof of Proposition \ref{prop:cvg of EM}.

\end{document}